\documentclass[12pt,reqno,a4paper]{amsart}
\usepackage{amscd,amsfonts,amssymb,amsmath,amsthm}
\usepackage[arrow,matrix]{xy}
\usepackage{graphicx,epstopdf}
\usepackage{subfigure,tikz}
\usepackage{color}
\usepackage{mathrsfs,euscript,hyperref}

\def\t{\theta}
\newtheorem{theorem}{Theorem}[section]
\newtheorem{lemma}[theorem]{Lemma}
\newtheorem{proposition}[theorem]{Proposition}

\theoremstyle{remark}
\newtheorem{remark}{Remark}[section]

\newtheorem{example}{Example}[section]
\theoremstyle{definition}
\newtheorem{definition}{Definition}[section]

\newcommand{\myp}{\mbox{$\:\!$}}
\newcommand{\mypp}{\mbox{$\;\!$}}

\newcommand{\myn}{\mbox{$\;\!\!$}}
\newcommand{\mynn}{\mbox{$\:\!\!$}}

\newcommand{\rme}{\mathrm{e}}

\newcommand{\xbar}[1]{%
   \kern0.15ex\hbox{%
     \vbox{%
       \hrule height 0.4pt 
       \kern0.3ex
       \hbox{%
         \kern-0.15em
         \ensuremath{#1}%
         \kern-0.1em
       }%
     }%
   \kern0.25ex}%
}

\newcommand{\xbarscr}[1]{%
   \kern0.2ex\hbox{%
     \vbox{%
       \hrule height 0.3pt 
       \kern0.2ex
       \hbox{%
         \kern-0.1em
         \ensuremath{{}_{#1}}%
         \kern-0.15em
       }%
     }%
   \kern0.25ex}%
}

\newcommand{\xbarscrscr}[1]{%
   \kern0.2ex\hbox{%
     \vbox{%
       \hrule height 0.3pt 
       \kern0.2ex
       \hbox{%
         \kern-0.1em
         \ensuremath{{}_{{}_{#1}}}%
         \kern-0.15em
       }%
     }%
   \kern0.25ex}%
}

\newcommand{\checkbxi}{\lefteqn{\boldsymbol{\xi}}\kern.13pc\check{\phantom{\xi}}\kern-.09pc}

\numberwithin{equation}{section}

\usepackage{cancel}
\usepackage[normalem]{ulem}

\begin{document}

\title[Gibbs measures of Potts model]{Gibbs measures of Potts model on Cayley trees: a survey and applications}

\author{U.A. Rozikov}

 \address{U.\ A.\ Rozikov \begin{itemize}
 \item[] V.I.Romanovskiy Institute of Mathematics of Uzbek Academy of Sciences;
\item[] AKFA University, 1st Deadlock 10, Kukcha Darvoza, 100095, Tashkent, Uzbekistan;
\item[] Faculty of Mathematics, National University of Uzbekistan.
\end{itemize}} \email {rozikovu@yandex.ru}

\begin{abstract} In this paper we give a systematic
review of the theory of Gibbs measures of Potts model  on Cayley trees (developed since 2013)  and discuss many applications of the Potts model
to real world situations: mainly biology, physics,  and some examples of
 alloy behavior, cell sorting, financial engineering, flocking birds, flowing foams,
image segmentation, medicine, sociology etc.
 \end{abstract}

 \maketitle

\section{Introduction}

The Potts model is defined by a Hamiltonian (energy) of configurations
of spins which take one of $q$ possible values on vertices of a lattice.
The model is used to
study the behavior of systems having multiple states (spins, colors, alleles, etc).
Since the model has a rich mathematical formulation it has been studied extensively.

Usually the results of investigations of a system (in particular, the Potts model) are presented
by a measure assigning a number to each (interesting) suitable properties of the system.
The Gibbs measure is one of such important measures in many problems of probability theory and statistical mechanics.
It is the measure associated with the Hamiltonian of a (biological, physical) system.
Each Gibbs measure gives a state of the system.

The main problem for a given Hamiltonian on a countable lattice is to describe its
all possible Gibbs measures. In case of uniqueness of the Gibbs measure
(for all values of the parameters, in particular a parameter can be temperature),
the system does not change its state. The existence of some values of parameters
at which the uniqueness of Gibbs measure switches to non-uniqueness is interpreted as phase transition
(the system will change its state).

The book \cite{R} presents all known (since 2012) mathematical results on (extreme) Gibbs measures on Cayley trees.
At that time, the Potts model on trees was not well studied, compared to the Ising model. Even results about
translation-invariant (simple) Gibbs measures were not complete and periodic Gibbs measures were not studied.

First complete result about translation-invariant Gibbs  measures of the Potts model
on Cayley trees appeared in 2014, \cite{KRK}.
Periodic and weakly periodic measures were studied after 2015.  These results mainly have applications to physics
as thermodynamical systems.

In this review, I have collected all these recent results and have presented many applications in biology,
medicine, sociology, physics etc.

\section{Preliminaries}

{\bf The Cayley tree $\Gamma^k$:} The Cayley tree $\Gamma^k$
 of order $ k\geq 1 $ is an infinite tree,
i.e., a graph without cycles, such that exactly $k+1$ edges
originate from each vertex. Let $\Gamma^k=(V,L,i)$, where $V$ is the
set of vertices $\Gamma^k$, $L$ the set of edges and $i$ is the
incidence function setting each edge $l\in L$ into correspondence
with its endpoints $x, y \in V$. If $i (l) = \{ x, y \} $, then
the vertices $x$ and $y$ are called the {\it nearest neighbors},
denoted by $l = \langle x, y \rangle $.

Fix a vertex $x^0\myn\in V$, interpreted as the \emph{root} of
the tree. We say that $y\in V$ is a \emph{direct successor}  of $x\in
V$ if $x$ is the penultimate vertex on the unique path leading from
the root $x^0$ to the vertex $y$; that is,
$d(x^0,y)=d(x^0,x)+1$ and $d(x,y)=1$. The set of all direct
successors of $x\in V$ is denoted $S(x)$.

For a fixed $x^0\in V$ we set $ W_n = \ \{x\in V\ \ | \ \ d (x,
x^0) =n \}, $
\begin{equation}\label{lp*}
 V_n = \ \{x\in V\ \ | \ \ d (x, x^0) \leq n \},\ \ L_n = \ \{l =
\langle x, y\rangle \in L \  | \ x, y \in V_n \}.
\end{equation}
For $x\in W_{n}$ the set $S(x)$ then has the form
\begin{equation}\label{sx}
S(x)=\{y\in W_{n+1}: \langle x, y\rangle\}.
\end{equation}

For any $x\in V$ denote
$$
W_m(x)=\{y\in V: d(x,y)=m\}, \ \ m\geq 1.
$$
Note that the sequence of balls $(V_n)$ ($n\in\mathbb{N}_0$)  is
\emph{cofinal} (see \cite[Section~1.2, page~17]{Ge}), that is, any
finite subset $\varLambda\subset V$ is contained in some $V_n$.

{\it Group representation of the tree.}

Let $G_k$ be a free product of $k + 1$ cyclic groups of the
second order with generators $a_1, a_2,\dots, a_{k+1}$,
respectively, i.e. $a_i^2=e$, where $e$ is the unit element.

It is known that there exists a one-to-one correspondence between the set of vertices $V$ of the
Cayley tree $\Gamma^k$ and the group $G_k$ (see Chapter 1 of \cite{R}
for properties of the group $G_k$).

Consideration of the Cayley trees has many
motivations (see e.g. \cite{Br}, \cite{DMR}, \cite{Dembo2}, \cite{Galanis}, \cite{Mos}, \cite{Rp}, \cite{Rr} and references therein).\\

{\bf Configuration space:}\label{cs}
Take a finite set $\Phi=\{1,2,...,q\}$, $q\geq 2$. This is the set of spin values.
 For $A\subseteq V$ a spin {\it configuration} $\sigma_A$
 on $A$ is defined as a function $$x\in
A\to\sigma_A(x)\in\Phi.$$

The set of all configurations coincides with $\Omega_A=\Phi^{A}$.
The cardinality of $\Omega_A$ is $|\Omega_A|=|\Phi|^{|A|}=q^{|A|}.$

We denote $\Omega=\Omega_V$ and
$\sigma=\sigma_V.$ Since $V$ is a countable set, $\Omega$ is an uncountable set.

Let $G^*$ be a
subgroup of the group $G_k$. A configuration $\sigma\in \Omega$ is
called $G^*$-periodic if $\sigma(yx)=\sigma(x)$ for any $x\in G_k$ and
$y\in G^*.$

  A configuration that is invariant with respect to all
shifts is called {\it translation-invariant}. \\

{\bf The $q$-state Potts model:} Consider the Potts model, where
the spin at each vertex $x\in V$ can
take values in the set $\Phi:=\{1,\dots,q\}$.

The Potts model with a \emph{nearest-neighbour interaction kernel}
$\{J_{xy}\}_{x,\myp y\myp\in\myp V}$ (i.e., such that
$J_{xy}=J_{yx}$ and $J_{xy}=0$ if $d(x,y)\ne1$) is defined by the
formal Hamiltonian
\begin{equation}\label{eq:H}
H(\sigma)=-\sum_{\langle x,\myp y\rangle\in L} \!J_{xy}
\mypp\delta_{\sigma(x),\myp\sigma(y)}-\sum_{x\in V}
\xi_{\sigma(x)}(x),\qquad \sigma\in\varPhi^V,
\end{equation}
where $\delta_{ij}$ is the Kronecker delta symbol, i.e.,
$$\delta_{ij}=\left\{\begin{array}{ll}
1, \ \ \mbox{if} \ \ i=j\\[2mm]
0, \ \ \mbox{if} \ \ i\ne j
\end{array}\right.$$
and $\xi(x)=(\xi_{1}(x),\dots,\allowbreak
\xi_{q}(x))\allowbreak \in \mathbb{R}^q$ is the external (possibly
random) field.

For $q=2$, the Potts model is equivalent to the Ising model (see Chapter 2 of \cite{R}).\\

{\bf Gibbs measure for the Potts model:} For each finite subset $\varLambda\subset V$
($\varLambda\ne\emptyset$) and any fixed sub-configuration
$\eta\in\varPhi^{\varLambda^c}$ (called the \emph{configurational
boundary condition}), the \emph{Gibbs measure} \index{Gibbs measure}
$\gamma^\eta_{\myn\varLambda}$ is a probability measure in
$\varPhi^\varLambda$ defined by the formula
\begin{equation}\label{eq:GD}
\gamma^\eta_{\myn\varLambda}(\varsigma)=\frac{1}{Z^\eta_{\myn\varLambda}(\beta)}\exp\biggl\{-\beta
H_{\myn\varLambda}(\varsigma)+\beta\sum_{x\in\varLambda}\sum_{y\in
\varLambda_{\vphantom{t}}^{c}}\!J_{xy}\mypp\delta_{\varsigma(x),\myp\eta(y)}\biggr\},\qquad
\varsigma\in\varPhi^{\varLambda},
\end{equation}
where  $H_{\myn\varLambda}$ is the restriction
of the Hamiltonian (\ref{eq:H}) to configurations in
$\varLambda$,
\begin{equation}\label{eq:H-Lambda}
H_{\myn\varLambda}(\varsigma)=-\sum_{\langle x,\myp y\rangle\in
L_\varLambda}\!\!J_{xy}
\mypp\delta_{\varsigma(x),\myp\varsigma(y)}-\sum_{x\in \varLambda}
\xi_{\varsigma(x)}(x),\qquad \varsigma\in\varPhi^{\varLambda},
\end{equation}
and $Z^\eta_{\myn\varLambda}(\beta)$ is the partition function
(or normalizing constant) \index{partition function}
$$
Z^\eta_{\myn\varLambda}(\beta)=\sum_{\varsigma\in\varPhi^{\varLambda}}
\exp\myn\Biggl\{-\beta
H_{\myn\varLambda}(\varsigma)+\beta\sum_{x\in\varLambda}\sum_{y\in
\varLambda_{\vphantom{t}}^{c}}J_{xy}\mypp\delta_{\varsigma(x),\myp\eta(y)}\Biggr\}.
$$

Now a measure $\mu=\mu_{\beta,\myp\xi}$ on $\varPhi^V$ is
called a \emph{Gibbs measure} if, \index{Gibbs measure}  for any non-empty finite set
$\varLambda\subset V$ and any $\eta\in\varPhi^{\varLambda^c}\!$,
\begin{equation}\label{eq:Gibbs}
\mu(\sigma_{\myn\varLambda}=\varsigma\mypp|\mypp\sigma_{\myn\varLambda^{c}}\myn=\eta)\equiv
\gamma_{\myn\varLambda}^{\eta}\myn(\varsigma),\qquad
\varsigma\in\varPhi^{\varLambda}.\\
\end{equation}

\textbf{The main problem:} is  to study the structure of the set $\mathcal G(H)$
of all Gibbs measures corresponding to a given Hamiltonian $H$.

The existence of Gibbs measures for a wide class of Hamiltonians was established in the
ground-breaking work of Dobrushin (see, e.g., \cite{Ge}, \cite{FV}, \cite{S}). However, a complete analysis
of the set of Gibbs measures for a given Hamiltonian is often a difficult problem.

A measure $\mu\in
\mathcal G(H)$ is called \emph{extreme}
(also called pure state) if it cannot be expressed as
$\mu=\lambda\myp \mu_1+(1-\lambda)\myp\mu_2$ for some
$\mu_1,\mu_2\in\mathcal G(H)$ with $\mu_1\ne\mu_2$.

The set of all
extreme measures in \strut{}$\mathcal G(H)$ \strut{}denoted by
${\rm ex}\mathcal G(H)$ is a \emph{Choquet simplex},  in the sense that any
$\mu\in\mathcal G(H)$ can be represented \strut{}as $\mu =
\int_{{\rm ex}\mathcal G(H)}\nu\,\rho({\mathrm d}{\nu})$, with some probability
measure $\rho$ on ${\rm ex}\mathcal G(H)$.

Thus the analysis of $\mathcal G(H)$ is reduced
to description of its extremal elements. Extremal Gibbs measures are very important to understanding all
possible local behaviors of the (biological and physical) system.

In this paper  we give some extreme
points of the set of Gibbs measures for the Potts model.\\

{\bf A method to describe Gibbs measure on trees:}
Follow \cite{BR} and \cite{KR} to explain the method of Markov random
field theory and its recurrent equations.

For a vector field
$V\ni x\mapsto {\bf h}(x)=(h_{1}(x),\dots, h_{q}(x))\in
\mathbb{R}^q$ and each $n\in\mathbb{N}_0=\{0,1,\dots\}$, define a probability
measure in $\varPhi^{V_n}$ by the formula
\begin{equation}\label{lp*}
\mu^h_n(\sigma_n)=\frac{1}{Z_n}\exp\left\{-\beta
H_n(\sigma_n)+\beta\!\sum_{x\in W_n}
\!h_{\sigma_n(x)}(x)\right\},\qquad \sigma_n\in\varPhi^{V_n},
\end{equation}
where $Z_n=Z_n(\beta, {\bf h})$ is the normalizing factor and
$H_n:=H_{V_n}$, that is (see~(\ref{eq:H-Lambda})),
\begin{equation}\label{eq:Hn}
H_n(\sigma_n)=-\sum_{\langle x,y\rangle\in L_n} \!\!J_{xy}
\mypp\delta_{\sigma_n(x),\myp\sigma_n(y)}-\sum_{x\in
V_n}\xi_{\sigma_n(x)}(x),\qquad \sigma_n\in\varPhi^{V_n}.
\end{equation}
The vector field $\{{\bf h}(x)\}_{x\in V}$ in (\ref{lp*}) is
called \emph{generalized boundary conditions (GBC)}.

One says that the probability distributions (\ref{lp*}) are
\emph{compatible} (and the intrinsic GBC $\{{\bf h}(x)\}$ are
\emph{permissible}) if for each $n\in\mathbb{N}_0$ the following
identity holds,
\begin{equation}\label{lp**}
\sum_{\omega\in
\varPhi^{W_{n+1}}}\!\mu^h_{n+1}(\sigma_n\myn\mynn\vee \omega)\equiv
\mu^h_{n}(\sigma),\qquad \sigma_n\in \varPhi^{V_{n}},
\end{equation}
where the symbol $\vee$ stands for concatenation of
sub-configurations.

By Kolmogorov's extension theorem (see, e.g., \cite{FV}, \cite[Chapter~II,
\S\mypp3, Theorem~4, page~167]{Shiryaev} and more suitable for our
setting is \cite[Theorem 6.2]{FV}), the compatibility
condition (\ref{lp**}) ensures that there exists a unique measure
$\mu^h=\mu^h_{\beta,\myp\xi}$ on $\varPhi^V$ such that, for all
$n\in\mathbb{N}_0$,
\begin{equation}\label{eq:mu-h}
\mu^h(\sigma_{V_n}\myn\mynn=\sigma_n)\equiv\mu^h_n(\sigma_n),\qquad
\sigma_n\in \varPhi^{V_n},
\end{equation}
or more explicitly (substituting~(\ref{lp*})),
\begin{equation}\label{eq:mu-h-ex}
\mu^h(\sigma_{V_n}\!=\sigma_n)=\frac{1}{Z_n}\exp\left\{-\beta
H_n(\sigma_n)+\beta\sum_{x\in W_n}h_{\sigma_n(x)}(x)\right\},\qquad
\sigma_n\in \varPhi^{V_n}.
\end{equation}
\begin{definition}
Measure $\mu^h$ satisfying \eqref{eq:mu-h} is called a
\emph{splitting Gibbs measure} (\emph{SGM}).
\end{definition}
\begin{remark}\label{rm:2.2}
Note that adding a constant $c=c(x)$ to all coordinates $h_i(x)$ of
the vector ${\bf h}(x)$ does not change the probability
measure (\ref{lp*}) due to the normalization $Z_n$. The same is true
for the external field ${\bf \xi}(x)$ in the Hamiltonian
(\ref{eq:Hn}). Therefore, without loss of generality consider
\emph{reduced GBC} ${\bf \check{h}}(x)$, for example defined \index{reduced GBC}
as
\begin{equation*}
\check{h}_i(x)=h_i(x)-h_q(x),\qquad i=1,\dots,q-1.
\end{equation*}
The same remark also applies to the external field
${\bf \xi}$ and its reduced version
${\bf \check{\xi}}(x)$, defined by
\begin{equation*}
\check{\xi}_i(x):=\xi_i(x)-\xi_q(x),\qquad i=1,\dots,q-1.
\end{equation*}
Such a reduction can equally be done by subtracting any
other coordinate,
$$
{}_{\ell}\check{h}_i(x):=h_i(x)-h_\ell(x),\qquad
{}_{\ell}\check{\xi}_i(x):=\xi_i(x)-\xi_\ell(x)\qquad (i\ne \ell).
$$
\end{remark}

Therefore, when working with vectors and
vector-valued functions and fields it will often be convenient to
pass from a generic vector $\boldsymbol{u}=(u_1,\dots,u_q)$ to a ``reduced vector''
$\boldsymbol{\check{u}}=(\check{u}_1,\dots,\check{u}_{q-1})\in
\mathbb{R}^{q-1}$ by setting $\check{u}_i:=u_i-u_q$
($i=1,\dots,q-1$).

The following statement describes a
criterion for the GBC $\{\boldsymbol{h}(x)\}_{x\in V}$ to
guarantee compatibility of the measures
$\{\mu^h_n\}_{n\in\mathbb{N}_0}$.

\begin{theorem}\label{ep} (see \cite{BR})
The probability distributions $\{\mu^h_n\}_{n\in\mathbb{N}_0}$
defined in \eqref{lp*} are compatible
 if and only if the
following vector identity holds
\begin{equation}\label{lp***}
\beta\myp\boldsymbol{\check{h}}(x)=\sum_{y\in
S(x)}\boldsymbol{F}\bigl(\beta\myp\boldsymbol{\check{h}}(y)+\beta\myp\boldsymbol{\check{\xi}}(y);
\myp\rme^{\beta J_{xy}}\bigr),\qquad x\in V,
\end{equation}
where
$\boldsymbol{\check{h}}(x)=(\check{h}_{i}(x),\dots,\check{h}_{q-1}(x))$,
\,$\boldsymbol{\check{\xi}}(x)=(\check{\xi}_{i}(x),\dots,\check{\xi}_{q-1}(x))$,
\begin{equation}\label{hxi}
\check{h}_{i}(x):=h_{i}(x)-h_{q}(x), \qquad \check
\xi_{i}(x):=\xi_{i}(x)-\xi_{q}(x),\qquad i=1,\dots,q-1,
\end{equation}
and the map
$\boldsymbol{F}(\boldsymbol{u};\theta)=(F_1(\boldsymbol{u};\theta),\dots,F_{q-1}(\boldsymbol{u};\theta))$
is defined for\\
$\boldsymbol{u}=(u_1,\dots,u_{q-1})\in\mathbb{R}^{q-1}$ and
$\theta>0$ by the formulas
\begin{equation}\label{eq:Fi}
F_i(\boldsymbol{u};\theta):=\ln\frac{(\theta-1)\mypp\rme^{u_i}+1+\sum_{j=1}^{q-1}\rme^{u_j}}{\theta+
\sum_{j=1}^{q-1}\rme^{u_j}},\qquad i=1,\dots, q-1.
\end{equation}
\end{theorem}
Now, using \cite[Theorem~(12.6)]{Ge} and according to the link between
GBC $\{\boldsymbol{h}(x)\}_{x\in V}$ and boundary laws (see \cite[Remark 1.4]{BR})
we make the following \textbf{crucial observations}:

\begin{itemize}
\item \ \  \emph{any extreme measure
$\mu\in{\rm ex}\mathcal G(H)$ is SGM}; therefore, the question of
uniqueness of the Gibbs measure is reduced to that in the SGM class. Moreover, for each given
temperature, the description
of the set $\mathcal G(H)$ is equivalent to the full description of the set of all extreme SGMs.
Therefore, in this paper we only interested to SGM on the Cayley trees.
\item \ \ \emph{Any SGM corresponds} to a solution of (\ref{lp***}) given in Theorem \ref{ep}.
Thus our main problem is reduced
to solving the functional equation (\ref{lp***}) and to check when a SGM corresponding to a solution is extreme.
\end{itemize}

To check the extremality of a Gibbs measure one can apply arguments used
for the reconstruction on trees \cite{FK}, \cite{Ke}, \cite{MarJ}, \cite{Mos2}, \cite{Mos}.

It is also known that a sufficient condition for non-extremality
(which is equivalent to solvability of the associated reconstruction) of a Gibbs measure is the Kesten-Stigum condition given in \cite{Ke}.

\section {Translation-invariant SGMs.}

In this section we consider the classic version of the Potts model:
\begin{equation}\label{ph}
H(\sigma)=-J\sum_{\langle x,y\rangle\in L}
\delta_{\sigma(x)\sigma(y)},
\end{equation}
where $J\in \mathbb R$, $\sigma(x)\in \Phi=\{1,\dots,q\}$ and
$\langle x,y\rangle$ stands for nearest neighbor vertices.

Let $S_q$ be the group of permutations on $\Phi$.

Take $\sigma\in \Phi^V$ and define an action of $\pi=(\pi_1, \dots, \pi_q)\in S_q$
on $\sigma$ (denoted by $\pi \sigma$) as
\[(\pi \sigma)(x)=\pi_{\sigma(x)}, \qquad \text{for all} \ \ x\in V.\]

Then it is easy to see that for any $\sigma\in \Phi^V$ and any $\pi \in S_q$ we have $H(\pi\sigma)=H(\sigma)$.

By Theorem \ref{ep} to each splitting Gibbs measure (SGM) of
the Hamiltonian (\ref{ph}) there is a vector-valued function ${h}_x$, such that
 for any $x\in V\setminus\{x^0\}$
the following equation holds:
\begin{equation}\label{p***}
 h_x=\sum_{y\in S(x)}F(h_y,\theta),
\end{equation}
where $F: h=(h_1, \dots,h_{q-1})\in \mathbb R^{q-1}\to F(h,\theta)=(F_1,\dots,F_{q-1})\in \mathbb R^{q-1}$ is defined as
$$F_i=\ln\left({(\theta-1)e^{h_i}+\sum_{j=1}^{q-1}e^{h_j}+1\over \theta+ \sum_{j=1}^{q-1}e^{h_j}}\right),$$
$\theta=\exp(J\beta)$, $S(x)$ is the set of direct successors of $x$ and $h_x=\left(h_{1,x},\dots,h_{q-1,x}\right)$ with
\begin{equation}\label{hh}
h_{i,x}={\tilde h}_{i,x}-{\tilde h}_{q,x}, \ \ i=1,\dots,q-1.
\end{equation}

Moreover, for any $h=\{h_x,\ \ x\in V\}$
satisfying (\ref{p***}) there exists a unique SGM $\mu$ for the Potts model.

In this section, we review main results of \cite{KRK}. Consider SGMs which are translation-invariant,
 i.e., assume $h_x=h=(h_1,\dots,h_{q-1})\in \mathbb R^{q-1}$ for all $x\in V$.
 Then from equation (\ref{p***}) we get $h=kF(h,\theta)$, i.e.,
\begin{equation}\label{pt}
h_i=k\ln\left({(\theta-1)e^{h_i}+\sum_{j=1}^{q-1}e^{h_j}+1\over \theta+ \sum_{j=1}^{q-1}e^{h_j}}\right),\ \ i=1,\dots,q-1.
\end{equation}
Denoting $z_i=\exp(h_i), i=1,\dots,q-1$,  from (\ref{pt}) we get
\begin{equation}\label{pt1}
z_i=\left({(\theta-1)z_i+\sum_{j=1}^{q-1}z_j+1\over \theta+ \sum_{j=1}^{q-1}z_j}\right)^k,\ \ i=1,\dots,q-1.
\end{equation}

\begin{remark}\label{rp} The permutation symmetry of Hamiltonian (\ref{ph}) mentioned above
consequents such symmetry to solutions of (\ref{pt1}): if $z=(z_1, \dots, z_{q-1})$ is a solution to
(\ref{pt1}) then by any permutation of its coordinates we get a solution of (\ref{pt1}) too.
\end{remark}

It is known that if $J<0$ (i.e. $\theta<1$) then for any $k\geq 1$, $q\geq 2$ the
anti-ferromagnetic Potts model has a unique TISGM (i.e., \cite[p.109]{R}, \cite{Ro9}) i.e.
the following theorem holds.

\begin{theorem}\label{antiT}
For the $q$-state anti-ferromagnetic ($J<0$) Potts model on the Cayley tree of order $k\geq 2$ there is unique TISGM.
\end{theorem}

The following theorem characterizes all solutions of (\ref{pt1}).
\begin{theorem}\label{tti} For any solution $z=(z_1,\dots,z_{q-1})$ of the system
of equations (\ref{pt1}) there exists $M\subset \{1,\dots,q-1\}$ and $z^*>0$ such that
$$z_i=\left\{\begin{array}{ll}
1, \ \ \mbox{if} \ \ i\notin M\\[3mm]
z^*, \ \ \mbox{if} \ \ i\in M.
\end{array}
\right.
$$
\end{theorem}

As a consequence of this theorem one can see that any TISGM of the Potts model corresponds to a solution of the following equation
\begin{equation}\label{rm}
z=f_m(z)\equiv \left({(\theta+m-1)z+q-m\over mz+q-m-1+\theta}\right)^k,
\end{equation}
for some $m=1,\dots,q-1$.

Put
$$T_{cr}={J\over \ln\left(1+{q\over k-1}\right)}.$$
In \cite{KRK} by analysis of solutions of the equation (\ref{rm}) the following theorem is proved.
\begin{theorem} \label{Theorem4}
For the $q$-state ferromagnetic ($J>0$) Potts model on the Cayley tree of
order $k\geq 2$ there are critical temperatures $T_{c,m}\equiv T_{c,m}(k,q)$, $m=1,\dots,[q/2]$
such that the following statements hold

\begin{itemize}
\item[1.] $$T_{c,1}>T_{c,2}>\dots>T_{c,[{q\over 2}]-1}>T_{c,[{q\over 2}]}\geq T_{cr};$$

\item[2.]
If $T>T_{c,1}$ then there exists a unique TISGM;

\item[3.]
If $T_{c,m+1}<T<T_{c,m}$ for some $m=1,\dots,[{q\over 2}]-1$ then there are
$1+2\sum_{s=1}^m{q\choose s}$ TISGMs.

\item[4.] If $T_{cr}\ne T<T_{c,[{q\over 2}]}$ then there are $2^{q}-1$ TISGMs.

\item[5.] If $T=T_{cr}$ then
the number of TISGMs is as follows
$$\left\{\begin{array}{ll}
2^{q-1}, \ \ \mbox{if} \ \ q -odd\\[2mm]
2^{q-1}-{q-1\choose q/2}, \ \ \mbox{if} \ \ q -even.
\end{array}\right.;$$

\item[6.] If $T=T_{c,m}$, $m=1,\dots,[{q\over 2}]$, \, ($T_{c,[q/2]}\ne T_{cr}$) then  there are
$1+{q\choose m}+2\sum_{s=1}^{m-1}{q\choose s}$ TISGMs.
\end{itemize}
\end{theorem}

\begin{figure}
\vspace{-.5pc} \centering
\includegraphics[width=13cm]{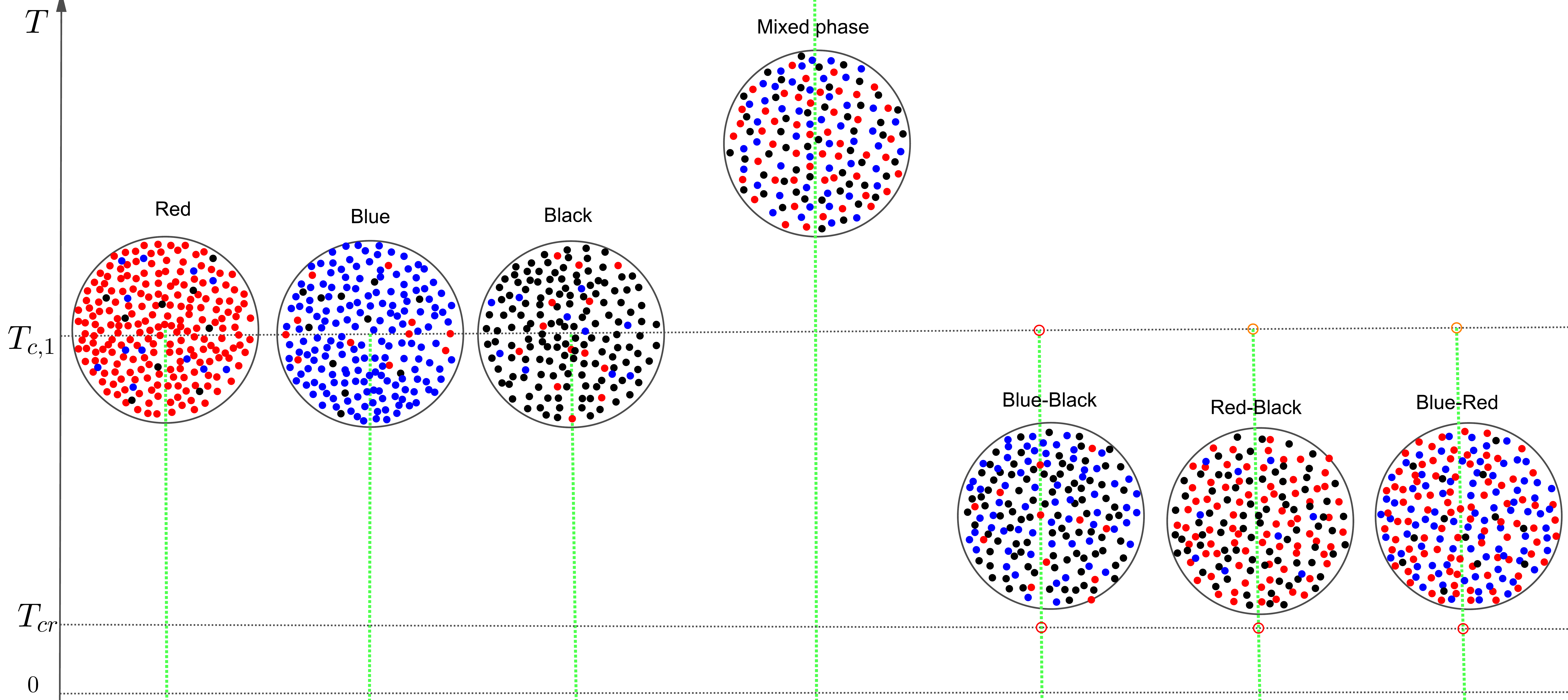}
\caption{3-state Potts model's typical configurations depending on temperature.  Here, states are: red=1, blue=2, black=3. All possible (seven)
cases are shown: In the mixed phase all three colors are seen with equal probability (with respect to the unique TISGM, i.e. free Gibbs measure).
The ``Red" phase (having many reds and very few other colors) is the typical configuration for the TISGM corresponding to the solution $(z_2,1,1)$.
Similarly, ``Blue" and ``Black" phases are the typical configurations for the TISGMs corresponding to the coordinate-permutations of $(z_2,1,1)$.
The ``Blue-Black" phase (having blue and black colors with equal probability, but very few red color) corresponds to $(z_1,1,1)$ and the remaining two phases correspond to its coordinate-permutations. Green vertical lines show the region of the temperature where the corresponding TISGM exists.
At the red points (shown at critical temperatures) on the blue lines the TISGMs do not exist.}\label{q3k2}
\end{figure}

\begin{remark} Let us make some useful remarks
\begin{itemize}
\item[1.] \ \ The critical temperature $T_{cr}={J\over \ln\left(1+{q\over k-1}\right)}$ is explicit for
 any $k\geq 2$ and any $q\geq 2$. For other critical values we know only existence and explicit formulas for $k=2$:
\begin{equation}\label{T2}
T_{c,m}\equiv T_{c,m}(2,q)={J\over \ln\left(1+2\sqrt{m(q-m)}\right)}, \ \ m=1,2,\dots [{q\over 2}].
\end{equation}
\item[2.] \ \ For the case $k=2$, $q=3$ we have two critical temperatures
$$T_{cr}={J\over \ln 4}, \ \ T_{c,1}={J\over \ln(1+2\sqrt{2})},$$
and up to $2^3-1=7$ TISGMs, denoted by $\mu_i$, $i=0,\dots,6$. Here, $\mu_0$ is free measure, corresponding to
solution 1.
Each TISGM defines its own measure of the cylinder $\{\sigma\in \Omega: \sigma(x)=j\}$, $x\in V$, $j=1,2,3$. Moreover,
 since these measures are translation invariant,
  the measure of this cylinder does not depend on $x\in V$:
 $$\mu_i(\{\sigma\in \Omega: \sigma(x)=j\})=p_{ij}.$$
Depending on the value of the solution $z$ of (\ref{rm}) one can calculate $p_{ij}$.
In case $z=1$, $p_{0j}=1/3$ (since $q=3$). If $\mu_i$ corresponds to a solution $z>1$ then for sufficiently low temperatures we have
$p_{i1}>2/3$ , and for $z<1$ we have $p_{i1}<1/3$.  Consequently, we can see
 typical configurations (of the seven TISGMs, i.e. phases) of the two-dimensional 3-state Potts model,
for different values of temperature as shown in Fig. \ref{q3k2}.
\end{itemize}
\end{remark}

\begin{theorem}\label{ext}  If $T\leq T_{c,1}$ then there are at least two extreme Gibbs measures for the $q$-state Potts model
on the Cayley tree of order $k\geq 2$.
\end{theorem}

Let ${\mathcal G}_{ti}$ be the set of all TISGMs.
The following gives relations between TISGMs.

\begin{theorem} Any measure $\mu\in {\mathcal G}_{ti}$ can not be non-trivial
convex mixture of measures from ${\mathcal G}_{ti}\setminus \{\mu\}$.
\end{theorem}

We now know that at sufficiently low
temperatures the maximal number of TISGMs is  $2^{q}-1$.
Moreover, it was shown that the number of TISGMs does not depend
on the order $k\geq 2$ of the Cayley tree.

 However it is not clear what kind of boundary
conditions (configurations)  are needed to get the
TISGMs as corresponding limits with the boundary conditions.

In \cite{GRR} we have answered this question.

\begin{remark}
1. It is known (see (\cite{Ge}, \cite{YH})) that the set $\mathcal G(J)$ of all
Gibbs measures (for a fixed parameter $J$ of the Hamiltonian (\ref{ph})) is a non-empty, compact convex set.
A limiting Gibbs measure is a Gibbs measure for the same $J$. Conversely, every extremal point of $\mathcal G(J)$
is a limiting Gibbs measure
with a suitable boundary condition for the same $J$.

2. It is known that any extreme Gibbs measure of a
Hamiltonian with nearest-neighbor interactions is a {\it splitting} Gibbs measure. Consequently,
any non-splitting Gibbs measure is not extreme. However, any splitting Gibbs measure
(not necessary extreme) is a limiting Gibbs measure, because  it corresponds
to a (generalized)\footnote{Recall that added a boundary field at each site of the
boundary is called a generalized boundary condition \cite{GRRR} or boundary law \cite{Ge}}
 boundary condition satisfying a compatibility (tree recursion) condition of Kolmogorov's theorem.

3. In \cite{C} it was shown that for non-extremal Gibbs measures on $\mathbb Z^d$
a Gibbs measure need not be a limiting Gibbs measure (see \cite{C} and \cite{FV} for more details).
\end{remark}

Consider the case $k=2$. In this case it is known that (see \cite{KRK}) for a given $m\leq [{q\over 2}]$,
there are vector solutions
$$(\underbrace{h_i,h_i,\dots,h_i}_m,\underbrace{0,0,\dots,0}_{q-m})$$
permuting coordinates of which one obtaines ${q\choose m}$ TISGMs.

Thus without loss of generality we can only consider
the measure $\mu_i(\theta,m)$ corresponding to vector
$${\bf h}(m,i)=(\underbrace{h_i,h_i,\dots,h_i}_m,\underbrace{0,0,\dots,0}_{q-m-1}).$$

Denote by $\mu_0\equiv\mu_0(\theta)$ the TISGM corresponding to solution $h_i\equiv 0$ and
by $\mu_i\equiv \mu_i(\theta,m)$ the TISGM corresponding to the solution $h_i(\theta, m)$,
$i=1,2$, $m=1,\dots, [{q\over 2}]$.

In \cite{GRR}  we have obtained all measures $\mu_i$ by changing boundary conditions (configurations).
To give the main result of \cite{GRR} we need the following result of \cite{KRK}.
Denote
\begin{equation}\label{tm}
\theta_m=1+2\sqrt{m(q-m)}, \ \ m=1,\dots,q-1.
\end{equation}
It is easy to see that
\begin{equation}\label{st}
\theta_m=\theta_{q-m} \ \ \mbox{and} \ \ \theta_{1}<\theta_2<\dots<\theta_{[{q\over 2}]-1}<\theta_{[{q\over 2}]}\leq q+1.
\end{equation}
Denote
\begin{equation}\label{s}
\begin{split}
x_1(m)={\theta-1-\sqrt{(\theta-1)^2-4m(q-m)}\over 2m},\\
x_2(m)={\theta-1+\sqrt{(\theta-1)^2-4m(q-m)}\over 2m}.
\end{split}
\end{equation}
\begin{proposition}\label{pw} Let $k=2$, $J>0$.
\begin{itemize}
\item[1.]
If $\theta<\theta_1$ then the system of equations (\ref{pt}) has
a unique solution $h_0=(0,0,\dots,0)$;

\item[2.]
If $\theta_{m}<\theta<\theta_{m+1}$ for some $m=1,\dots,[{q\over 2}]-1$ then the system of equations (\ref{pt}) has
solutions
$$h_0=(0,0,\dots,0), \ \ h_{1i}(s), \ \ h_{2i}(s), \ \ i=1,\dots, {q-1\choose s}, $$
$$h'_{1i}(q-s), \ \ h'_{2i}(q-s), \ \ i=1,\dots, {q-1\choose q-s}, \ \ s=1,2,\dots,m,$$
where $h_{ji}(s)$, (resp. $h'_{ji}(q-s)$)\, $j=1,2$ is a vector with $s$ (resp. $q-s$)
coordinates equal to $2\ln x_j(s)$ (resp. $2\ln x_j(q-s)$) and the remaining $q-s-1$ (resp. $s-1$)
coordinates equal to 0. The number of such solutions is equal to
$$1+2\sum_{s=1}^m{q\choose s};$$

\item[3.] If $\theta_{[{q\over 2}]}<\theta\ne q+1$ then there are $2^q-1$ solutions to (\ref{pt});

\item[4] If $\theta=q+1$ then the
number of solutions is as follows
$$\left\{\begin{array}{ll}
2^{q-1}, \ \ \mbox{if} \ \ q  \ \ \mbox{is odd}\\[2mm]
2^{q-1}-{q-1\choose q/2}, \ \ \mbox{if} \ \ q \ \ \mbox{is even};
\end{array}\right.$$

\item[5.] If $\theta=\theta_m$, $m=1,\dots,[{q\over 2}]$, \,($\theta_{[{q\over 2}]}\ne q+1$) then  $h_{1i}(m)=h_{2i}(m)$. The number of solutions is equal to
$$1+{q\choose m}+2\sum_{s=1}^{m-1}{q\choose s}.$$
\end{itemize}
\end{proposition}

\begin{remark}\label{sev} We note that
 \begin{itemize}
\item[1)] \ \ By Proposition \ref{pw} for $k=2$ and $J>0$ we have the
 {\rm full} description of solutions to the system of equations (\ref{pt}).
 Consequently, this gives the full description of TISGMs. Moreover,
 depending on parameter
 $\theta$ the maximal number of such measures can be $2^{q}-1$.

\item[2)] \ \ Recall $\theta_c={k+q-1\over k-1}$, $k\geq 2$, $q\geq 2$.
For $k=2$ by (\ref{tm}) we have
$$\left\{\begin{array}{ll}
\theta_c> \theta_m, \ \ \mbox{for all} \ \ m\in \{1,\dots,q-1\}\setminus \{q/2\}\\[2mm]
\theta_c=\theta_m, \ \ \mbox{for} \ \ m=q/2.
\end{array}
\right.$$
\end{itemize}
\end{remark}

Let $\omega\in \Omega$ be a configuration such that
$$c^l(\omega)=\sum_{s\in S(t)}\delta_{l\omega(s)}$$ is independent of $t\in V\setminus\{x^0\}.$

For a given $m\in\{1,\dots,[{q\over 2}]\}$ and $J>0$ introduce the following
sets of configurations:
$${\mathbb B}_{m}=\{\omega\in \Omega: c^1(\omega)=\dots =c^m(\omega),\ \ c^{m+1}(\omega)=\dots =c^{q-1}(\omega)=c^q(\omega)\},$$
$${\mathbb B}^+_{m,0}=\{\omega\in \mathbb B_m: c^1(\omega)>c^q(\omega)\},$$
$${\mathbb B}^0_{m,0}=\{\omega\in \mathbb B_m: c^1(\omega)=c^q(\omega)\},$$
$${\mathbb B}^-_{m,0}=\{\omega\in \mathbb B_m: c^1(\omega)<c^q(\omega)\},$$
 $${\mathbb B}^+_{m,1}=\{\omega\in \mathbb B_m: J\left(c^1(\omega)-c^q(\omega)\right)>h_1\},$$
  $${\mathbb B}^0_{m,1}=\{\omega\in \mathbb B_m: J\left(c^1(\omega)-c^q(\omega)\right)=h_1\},$$
$${\mathbb B}^-_{m,1}=\{\omega\in \mathbb B_m: J\left(c^1(\omega)-c^q(\omega)\right)<h_1\}.$$

Denote by $P^\omega$ the limiting Gibbs measure corresponding to a boundary configuration $\omega$.
The main result of \cite{GRR} is the following
\begin{theorem}\label{t2} The following assertions hold

1) If $\theta=\theta_m$, for some $m=1,\dots,[{q\over 2}]$ then
\begin{equation}\label{1a}
P^\omega=\left\{
\begin{array}{ll}
    \mu_1(\theta, m), & \hbox{\textit{if}}\ \ \omega\in \mathbb B_{m,1}^+\cup \mathbb B_{m,1}^0\\[2mm]
    \mu_0(\theta), & \hbox{\textit{if}}\ \ \omega\in \mathbb B^-_{m,1}
    \end{array}
\right.
\end{equation}

2) If $\theta_m<\theta<\theta_c=q+1$ then

\begin{equation}\label{1b}
P^\omega=\left\{
\begin{array}{lll}
    \mu_2(\theta, m), & \hbox{\textit{if}}\ \ \omega\in \mathbb B_{m,1}^+\\[2mm]
    \mu_1(\theta,m), & \hbox{\textit{if}}\ \ \omega\in \mathbb B^0_{m,1}\\[2mm]
     \mu_0(\theta), & \hbox{\textit{if}}\ \ \omega\in \mathbb B_{m,1}^-
\end{array}
\right.
\end{equation}

3) If $\theta=\theta_c$ then

\begin{equation}\label{1c}
P^\omega=\left\{
\begin{array}{llll}
    \mu_2(\theta, m), & \hbox{\textit{if}}\ \ \omega\in \mathbb B^{+}_{m,0}\\[2mm]
    \mu_0(\theta), & \hbox{\textit{if}}\ \ \omega\in \mathbb B^{-}_{m,0}\cup \mathbb B^0_{m,0}
\end{array}
\right.
\end{equation}

4) If $\theta>\theta_c$ then
\begin{equation}\label{1d}
P^\omega=\left\{
\begin{array}{lll}
    \mu_2(\theta, m), & \hbox{\textit{if}}\ \ \omega\in \mathbb B_{m,0}^+\\[2mm]
    \mu_1(\theta,m), & \hbox{\textit{if}}\ \ \omega\in \mathbb B^-_{m,0}\\[2mm]
     \mu_0(\theta), & \hbox{\textit{if}}\ \ \omega\in \mathbb B_{m,0}^0
\end{array}
\right. 
\end{equation}
\end{theorem}
For concrete examples of the boundary conditions mentioned in Theorem \ref{t2} see \cite{GRR}.

\section{Conditions of (non-) extremality of TISGMs}

In this section we review results of papers \cite{KR}, \cite{RKK} about sufficient conditions for (non-)extremality
of TISGMs of the Potts model,
depending on coupling strength parameterized by $\theta$, the block size $m$ and the branch of the boundary
law $z$.

Recall that by $\mu_i(\theta,m)$ we denote the TISGM which corresponds to
the values of $z_i(\theta,m)$, $i=1,2$, which is a solution to (\ref{rm}).\\

{\bf Non-extremality.}

Define the following numbers:
\begin{equation}\label{t*}
\widehat\theta=(\sqrt{2}-1)q+2m+1, \ \ \theta^*=1+(\sqrt{2}+1)q-2m.
\end{equation}

\begin{figure}
\includegraphics[width=9cm]{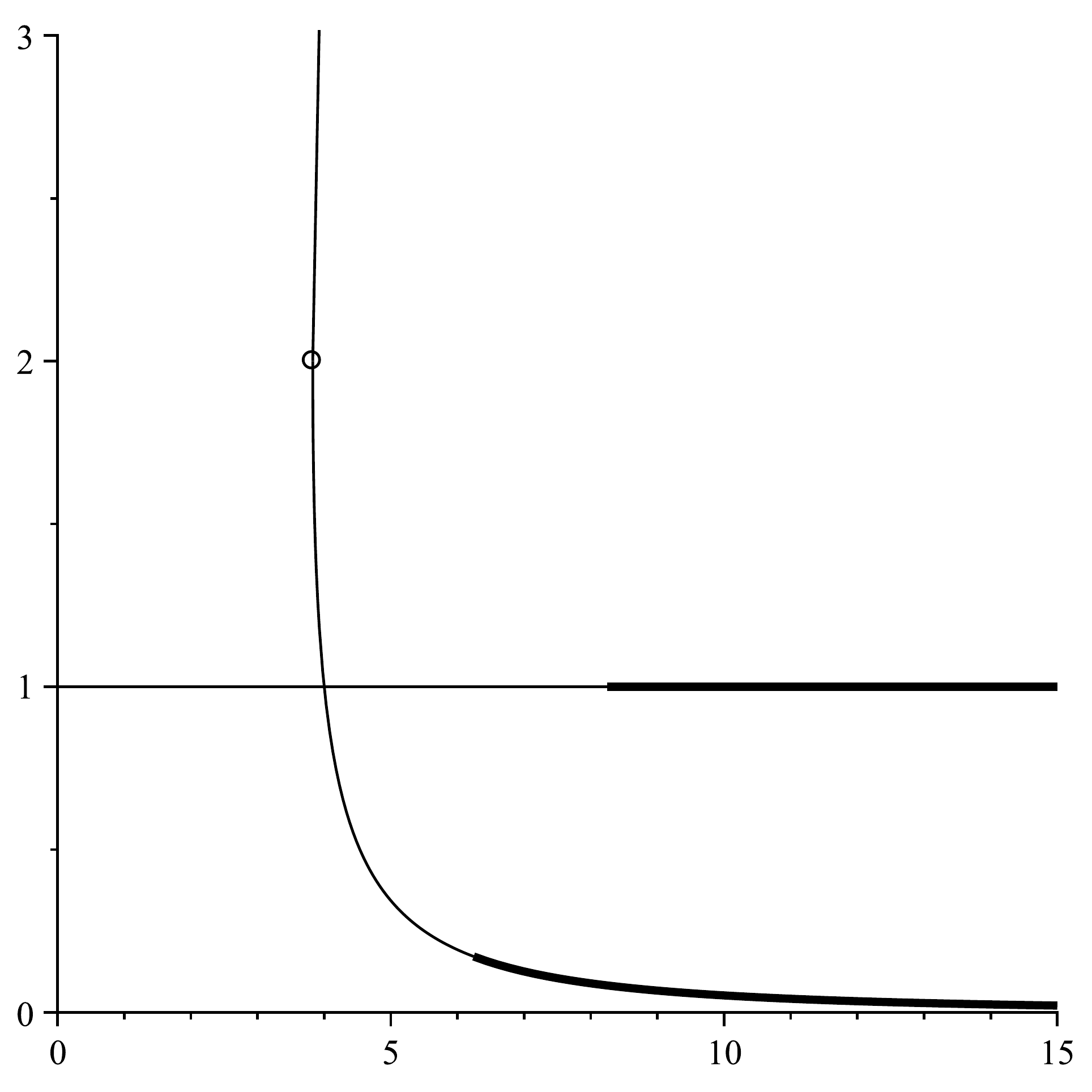}
\caption{ The graphs of the functions $z_i=z_i(m,\theta)$, $i=1,2$, for $q=3$ and $m=1$. The circle dot having coordinate $(\theta_1, 2)$ separates graph of $z_1$ from graph of $z_2$. The bold curves correspond to regions of solutions where the corresponding TISGM is non-extreme. This figure corresponds to Part (i), b) of Theorem \ref{tne}.}\label{fn2}
\includegraphics[width=9cm]{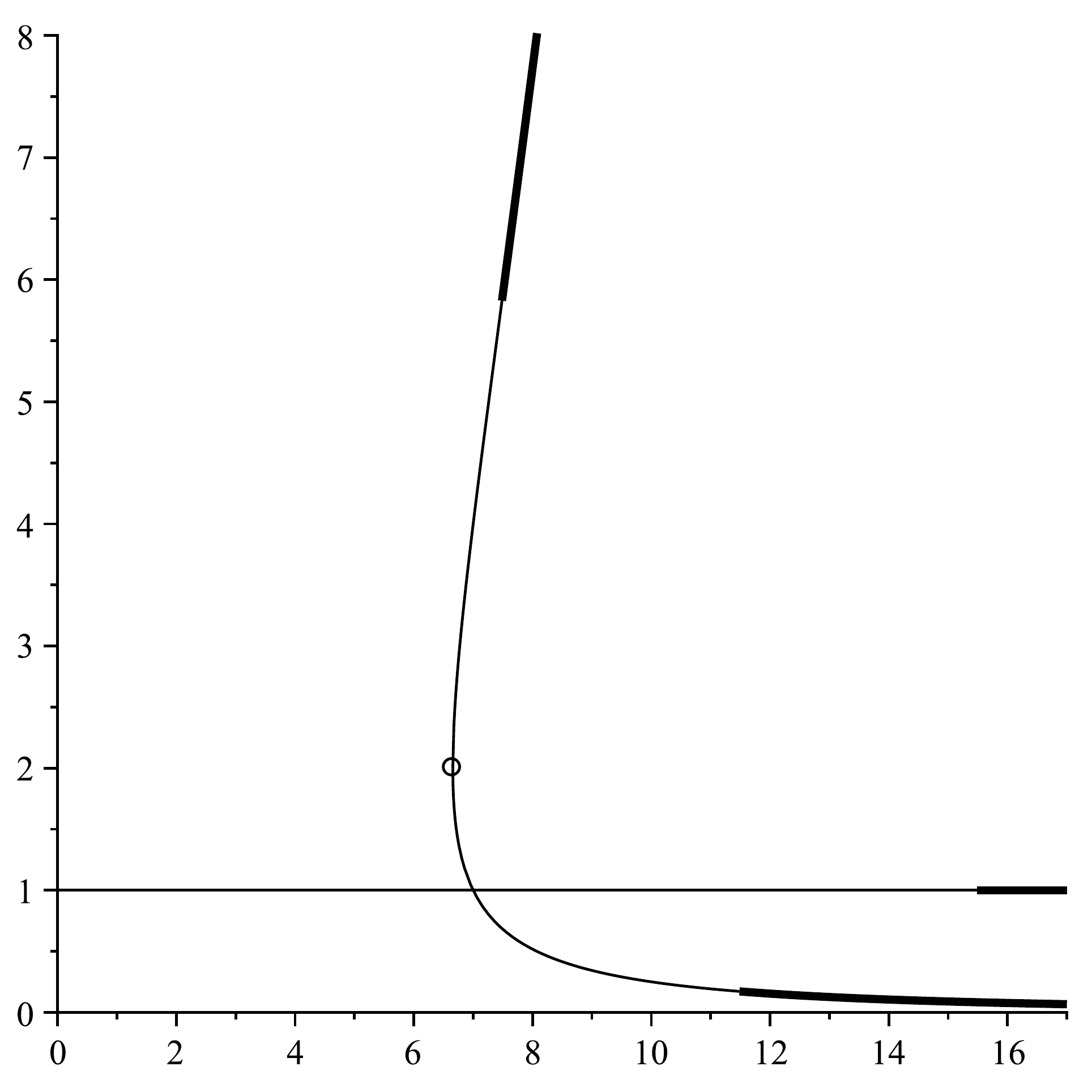}
\caption{The graphs of the functions $z_i=z_i(m,\theta)$, $i=1,2$, for $q=6$ and $m=2$. The circle dot having coordinate $(\theta_2, 2)$ separates graph of $z_1$ from graph of $z_2$. The bold lines correspond to regions of solutions where corresponding TISGM is non-extreme. The upper bold curve corresponds to part (ii), d) and the lower bold curve corresponds to part (i), b) of Theorem \ref{tne}.}\label{fn3}
\end{figure}

\begin{figure}
\includegraphics[width=9cm]{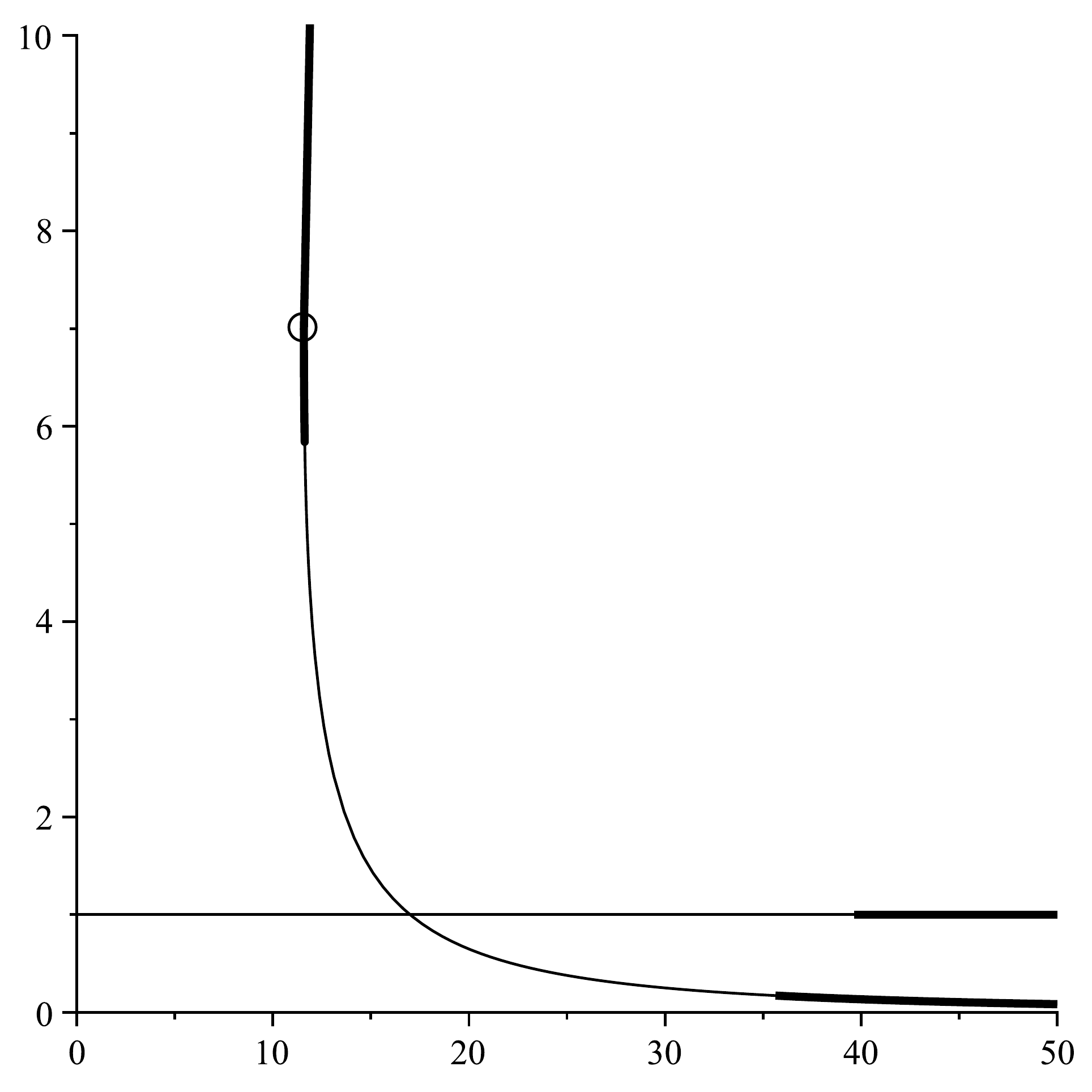}
\caption{The graphs of the functions $z_i=z_i(m,\theta)$, $i=1,2$, for $q=16$ and $m=2$. The circle dot having coordinate $(\theta_2, 7)$ separates graph of $z_1$ from graph of $z_2$. The bold lines correspond to regions of solutions where corresponding TISGM is non-extreme. The bold $z_2$ corresponds to Part (ii), c); Bold parts of $z_1$ correspond to Part (i), a) (upper bold curve) and Part (i), b) (lower bold curve) of Theorem \ref{tne}.}\label{fn4}
\end{figure}

\begin{theorem}\label{tne} Let $k=2$, $2m<q$. Then the following statements hold.
\begin{itemize}
\item[(i)] \ \ \ Assume one of the following conditions is satisfied:

\begin{itemize}
\item[a)] $2\leq m\leq q/7$ and $\theta\in [\theta_m, \widehat\theta)$;

\item[b)] $\theta\in (\theta^*,+\infty)$.
\end{itemize}

Then  $\mu_1(\theta,m)$ is non-extreme.

\item[(ii)] \ \ \ Assume one of the following conditions is satisfied:

\begin{itemize}
\item[c)] $2\leq m\leq q/7$ and $\theta\geq \theta_m$;

\item[d)]  $q<7m$, $m\geq 2$ and $\theta\in (\widehat\theta, +\infty)$.

\end{itemize}

Then  $\mu_2(\theta,m)$ is non-extreme.(See Fig.\ref{fn2}-\ref{fn4})\\
\end{itemize}
\end{theorem}

{\bf Extremality.}
Denote
$\theta_1=1+2\sqrt{q-1}$ and $\theta^*=1+(\sqrt{2}+1)q-2m$.

\begin{figure}
\includegraphics[width=8cm]{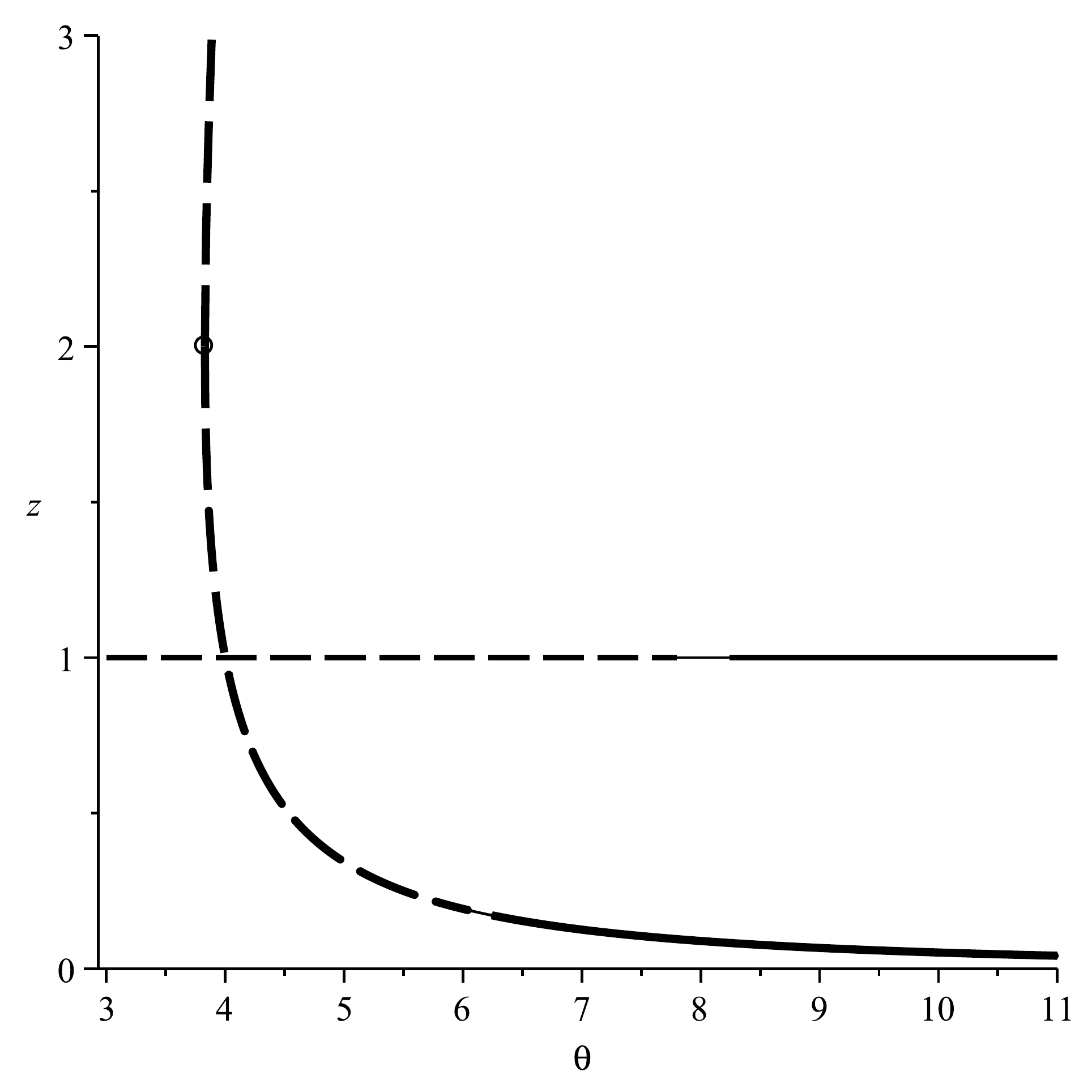}
\caption{The graphs of the functions $z_i=z_i(m,\theta)$, for $q=3$, $m=1$ and the graph of $z(\theta)\equiv 1$. The bold curves correspond to regions of solutions where the corresponding TISGM is non-extreme (corresponds to part (ii), b) of Theorem \ref{tne}). The dashed bold curves correspond to regions of solutions where the corresponding TISGM is extreme (see parts (a) and (b) of Theorem \ref{t1}). The gaps between the two types of curves are given by thin curves.}\label{fn5}
\end{figure}

\begin{figure}
\includegraphics[width=8cm]{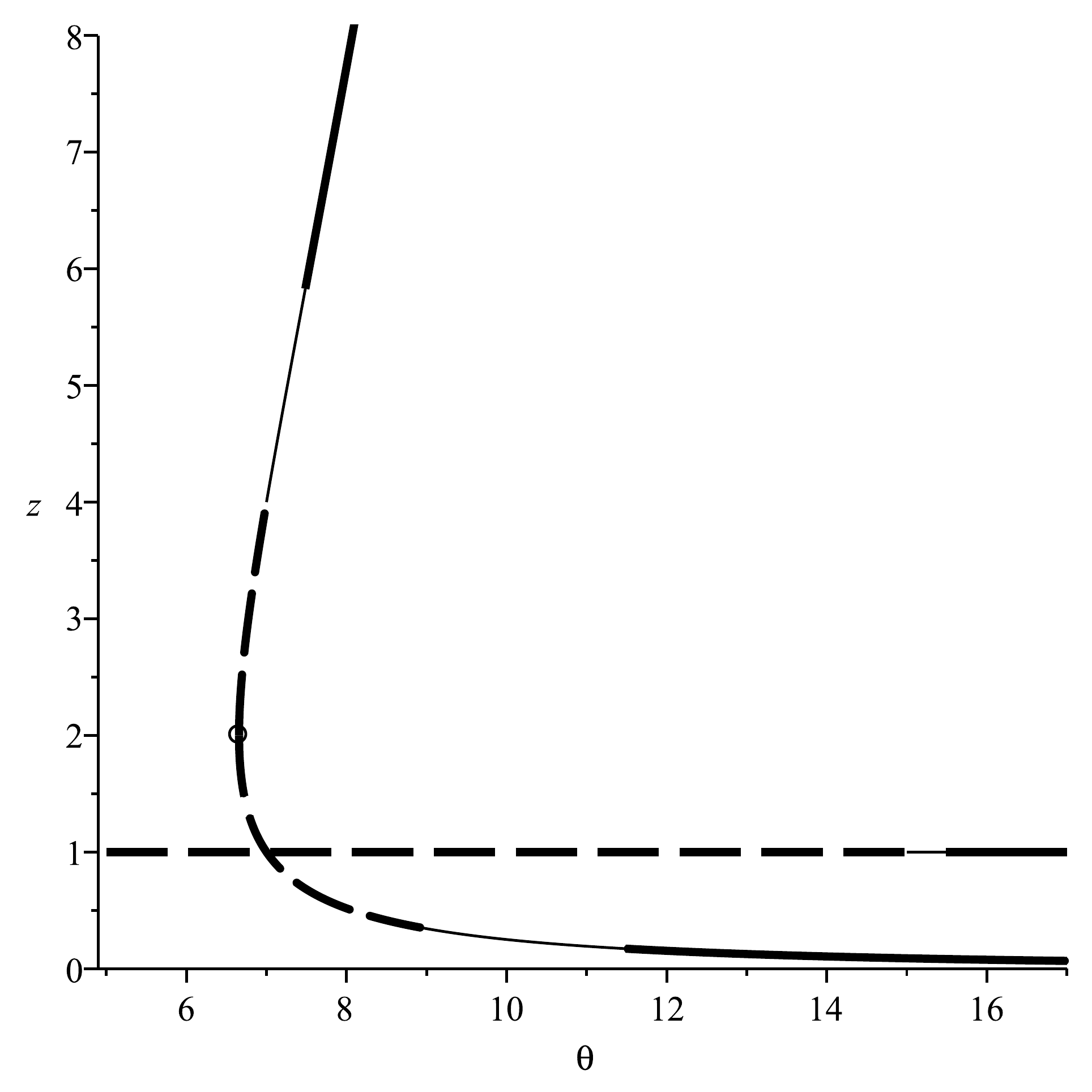}
\caption{The graphs of the functions $z_i=z_i(m,\theta)$, for $q=6$, $m=2$ and the graph of $z(\theta)\equiv 1$. The types of curves corresponding
to certain extremality and certain non-extremality are as in Fig.\ref{fn5} (corresponds to Theorem \ref{t3}).}\label{fn7}
\end{figure}

\begin{theorem}\label{t1} If $k=2$, $m=1$ then the following is true.
\begin{itemize}
\item[(a)] 
\ \ \ - If $q=3,4,\dots,16$ then there exists $\theta^{**}$ such that $\theta_c=q+1<\theta^{**}<\theta^*$ and the
measure $\mu_1(\theta,1)$ is extreme for any $\theta\in [\theta_1, \theta^{**})$.
Moreover $\theta^{**}$ is the unique positive solution of the following equation
$$\theta^3-(q-3)\theta^2-(2q-7)\theta-(q+5)=0.$$
 - If $q\geq 17$ then there are $\bar\theta_1, \bar\theta_2\in (\theta_1,\theta_c)$ such that $\bar\theta_1<\bar\theta_2$ and
the measure $\mu_1(\theta,1)$ is extreme for any $\theta\in [\theta_1, \bar\theta_1)\cup (\bar\theta_2,\theta^{**})$.
Moreover $\bar\theta_1, \bar\theta_2$ are positive solutions of the following equation
$$\theta^3-(q-1)\theta^2-(2q-3)\theta+(4q^2-13q+11)=0.$$
\item[(b)] \ \ \ The measure $\mu_2(\theta,1)$ is extreme for any $\theta\geq  \theta_1$, $q\geq 2$. (see Fig.\ref{fn5})
\end{itemize}
\end{theorem}
The following theorem corresponds to the case: $m\geq 2$. From condition $2\leq m\leq [q/2]$ it follows that $q\geq 4$.
\begin{theorem}\label{t3} Let $k=2$.

\begin{itemize}

\item[(i)] \ \ If $m=2$ then the following is true.

\item[(i.1)] \ \ \ If $m=2$ then for each $q=4,5,6,7,8$ there exists $\breve\theta>\theta_c=q+1$ such that the measure $\mu_1(\theta,2)$ is extreme for any $\theta\in [\theta_2, \breve\theta)$.
\item[(i.2)] \ \ \ For each $q\geq 9$ there exists $\theta^\dag\in (\theta_2,q+1)$ such that the measure $\mu_1(\theta,2)$ is extreme for any $\theta\in [\theta^\dag, \breve\theta)$, where $\theta^\dag=\theta^\dag(q)$ is the unique solution of
    $$ \theta^3-(q+3)\theta^2+(6q-17)\theta-(9q-19)=0$$ and $\breve\theta=\breve\theta(q)$ is the unique solution of
    $$\theta^3-(q+3)\theta^2-(2q-15)\theta-(q+13)=0.$$
\item[(ii)] \ \ \ If $m=2$ then for each $q=4,5,6,7,8$ there exists $\grave{\theta}=\grave{\theta}(q)$ such that $\theta_2<\grave{\theta}\leq q+1$ and $\mu_2(\theta,2)$ is extreme for $\theta\in [\theta_2, \grave\theta)$ (see Fig.\ref{fn7}).
\item[(iii)] \ \ \ If $q<{m+1\over 2m}\left[3m+1+\sqrt{m^2+6m+1}\right]$ and $m\geq 2$ then the measure $\mu_1(\theta_m,m)=\mu_2(\theta_m,m)$ is extreme.
\end{itemize}
\end{theorem}
Now following \cite{RKK} we give a theorem which is more general than Theorem \ref{t3}.

\begin{theorem}\label{rkk} Let $k=2, \ m\geq2$. Then
\begin{itemize}
\item[1.] \ \ If $2m\leq q<{m+1\over 2m}[3m+1+\sqrt{m^2+6m+1}]$, then there is
$\ddot{\theta}>q+1$ such that the measure $\mu_1(\theta,m)$
is extreme for any $\theta\in [\theta_m,\ddot{\theta})$;

\item[2.] \ \ If $q>{m+1\over 2m}[3m+1+\sqrt{m^2+6m+1}]$, then there exists
$\overline{\theta}\in(\theta_m,\theta_c)$ such that the measure
$\mu_1(\theta,m)$ is extreme for any $\theta\in
(\overline{\theta},\ddot{\theta})$;

\item[3.] \ \ If $2m\leq q<m+{1\over 4m}(m+1+ \sqrt{m^2+2m+7})^2$
there is  $\bar{\bar{\theta}}\in(\theta_m,+\infty)$ such that
the measure $\mu_2(\theta,m)$ is extreme for any $\theta\in
[\theta_m,\bar{\bar{\theta}})$.
\end{itemize}
\end{theorem}
\begin{remark} Note that
\begin{itemize}
\item \ \ If $q>m+{1\over 4m}(m+1+ \sqrt{m^2+2m+7})^2$, then
it is unknown whether the measure $\mu_2(\theta,m)$ is extreme in this case.

\item \ \ In \cite{RKK} for the 3-state Potts model on the Cayley tree of order $k=3$,
the extremality questions are studied for certain TISGMs. The case $k\geq 4$ is not considered yet. Because in this case
there is not explicit formula of nontrivial solutions to the equation (\ref{rm}).
\end{itemize}
\end{remark}

\begin{figure}
\includegraphics[width=8cm]{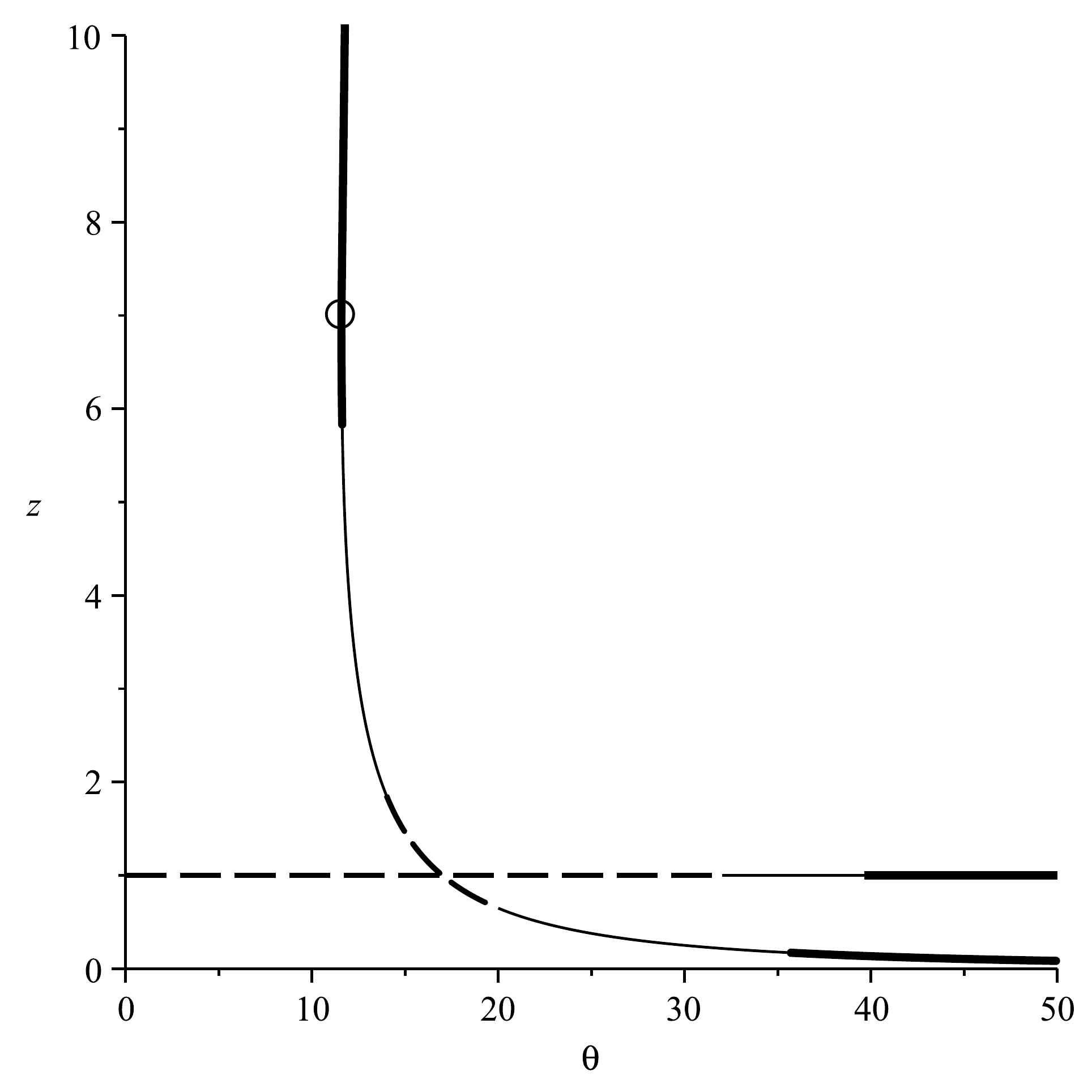}
\caption{The graphs of the functions $z_i=z_i(m,\theta)$, for $q=16$, $m=2$ and the graph of $z(\theta)\equiv 1$. The types of curves are again as in Fig.\ref{fn5} (corresponds to part a) of Theorem \ref{th2}).}\label{fn9}
\end{figure}

Now we give results which are valid for $\theta$ close to the critical value $\theta_c$
at which the lower branches of the boundary law $z$ degenerate into the free value $z=1$ and
the corresponding Markov chains become close to the free chain.
\begin{theorem}\label{th2} Let $k=2$.
 \begin{itemize}
 \item[a)] \ \ For each $m\leq [q/2]$ there exists a neighborhood $U_m(\theta_c)$ of $\theta_c$ such that the
measure $\mu_1(\theta,m)$ is extreme if $\theta\in U_m(\theta_c)$ (see Fig.\ref{fn9}).

\item[b)] \ \ If $m=1$ or condition (iii) of Theorem \ref{t3} is satisfied then there exists a neighborhood $V_m(\theta_m)$ of $\theta_m$ such that
measures $\mu_i(\theta,m),$ $i=1,2$ are extreme if $\theta\in V_m(\theta_m)$.
\end{itemize}
\end{theorem}
Recall that the model is called ferromagnetic if $\theta>1$, and as usual
a punctured neighborhood of $\theta_c$ is an open neighborhood from which the value
$\theta_c$ was removed.

\begin{theorem}\label{th8} For the ferromagnetic
$q$-state Potts model (with $q\geq 3$) on the Cayley tree
of order two, there exists a punctured neighborhood $U(\theta_c)$ of $\theta_c$
such that there are at least $2^{q-1}+q$ extreme TISGMs for each $\theta\in U(\theta_c)$.
\end{theorem}

\section{TISGMs in case of a non-zero external field}

This section is devoted to review results of \cite{BR} concerning TISGMs in case of a non-zero external field.


 For \strut{}simplicity,  in (\ref{eq:H}) we choose the external field such that
 all coordinates of the (reduced)
vector $\check{\boldsymbol{\xi}}{}^0\mynn$ (see Remark  \ref{rm:2.2}) are zero except one; due
to permutational symmetry, we may assume, without loss of
generality, that $\check{\xi}_1^0=\alpha\in\mathbb{R}$ and
$\check{\xi}^0_2=\dots=\check{\xi}^0_{q-1}=0$,
\begin{equation}\label{eq:xi0}
\boldsymbol{\check{\xi}}{}^0\mynn = (\alpha,0,\dots,0)\in
\mathbb{R}^{q-1}.
\end{equation}
We also write
$$
\boldsymbol{\check{h}}{}^0\mynn=(\check{h}_1^0,\dots,\check{h}_{q-1}^0)\in\mathbb{R}^{q-1}.
$$
Then, denoting $z_i:=\theta^{\mypp\check{h}_i^0/k}$
($i=1,\dots,q-1$), the compatibility equations \eqref{p***} are
equivalently rewritten in the form
\begin{equation}\label{pt1}
\left\{\begin{aligned} z_1&=1+\frac{(\theta-1)(\theta^{\alpha}
z_1^k-1)}{\theta+\theta^{\alpha}_{\vphantom{X}}
z_1^k+\sum_{j=2}^{q-1}z_j^k},\\
z_i&=1+\frac{(\theta-1)(z_i^k-1)}{\theta+\theta^{\alpha}_{\vphantom{X}}
z_1^k+\sum_{j=2}^{q-1}z_j^k},\qquad i=2,\dots,q-1.
\end{aligned}
\right.
\end{equation}

Assume $\alpha\ne 0$. The system \eqref{pt1} with $\theta>1$ is
reduced either to a single equation
\begin{equation}\label{rm=1-int}
u=1+\frac{(\theta-1)(\theta^\alpha
u^k-1)}{\theta+\theta^\alpha_{\vphantom{X}}
u^k+q-2}
\end{equation}
or to the system of equations (indexed by $m=1,\dots, q-2$)
\begin{equation}\label{uv-int}
\left\{\begin{aligned} u&=1+\frac{(\theta-1)(\theta^\alpha
u^k-1)}{\theta+\theta^\alpha_{\vphantom{X}}
u^k+m v^k+q-2-m},\\
1&=\frac{(\theta-1)(1+v+\dots+v^{k-1})}{\theta+\theta^\alpha_{\vphantom{X}}
u^k+m v^k+q-2-m},\\
v&\ne 1.
\end{aligned}
\right.
\end{equation}

To give the solvability of the
equation~\eqref{rm=1-int} let us introduce some notations.
Denote
\begin{equation}\label{eq:theta'c}
\theta_{\rm c}=\theta_{\rm c}(k,q):=
\frac{1}{2}\,\Biggl(\sqrt{(q-2)^2+4\myp(q-1)\left(\frac{k+1}{k-1}\right)^2}-(q-2)\Biggr).
\end{equation}
\begin{equation}\label{eq:b}
b= b(\theta):=\frac{\theta\myp(\theta+q-2)}{q-1}.
\end{equation}
 For $\theta\ge\theta_{\rm c}$, denote by $x_{\pm}=
x_{\pm}(\theta)$ the roots of the quadratic equation
\begin{equation}\label{eq:quadratic-intro}
(b+x)(1+x)=k\myp(b-1)\myp x
\end{equation}
with discriminant
\begin{equation}\label{eq:D}
D= D(\theta):=\bigl(k\myp(b-1)-(b+1)\bigr)^2-4b=
(b-1)(k-1)^2\left(b-\left(\frac{k+1}{k-1}\right)^2\right),
\end{equation}
that is,
\begin{equation}\label{eq:x+/--intro}
x_{\pm}= x_{\pm}(\theta):=\frac{(b-1)(k-1)-2\pm\sqrt{D}}{2}.
\end{equation}
Introduce the following notations
\begin{equation}
\label{eq:a_pm} a_{\pm}=
a_{\pm}(\theta):=\frac{1}{x_\pm}\left(\frac{1+x_\pm}{b+x_\pm}\right)^k,\qquad
\theta\ge\theta_{\rm c}.
\end{equation}
\begin{equation}\label{eq:alpha_pm}
\alpha_{\pm}=\alpha_{\pm}(\theta):=-(k+1)+\frac{1}{\ln\theta}\ln\frac{q-1}{a_{\mp}},\qquad
\theta\ge\theta_{\rm c},
\end{equation}
so that $\alpha_{-}(\theta_{\rm c})=\alpha_{+}(\theta_{\rm c})$ and
$\alpha_{-}(\theta)<\alpha_{+}(\theta)$ for $\theta>\theta_{\rm c}$.

\begin{theorem}\label{th:3.7}
Let $\nu_0(\theta,\alpha)$ denote the number of TISGMs corresponding to the solutions $u>0$ of
the equation~\eqref{rm=1-int}. Then
\begin{equation*}
\nu_0(\theta,\alpha)=\begin{cases} 1&\text{if}\ \, \theta\leq
\theta_{\rm c} \ \,\text{or} \ \,\theta>\theta_{\rm c} \ \,
\text{and} \ \, \alpha\notin [\alpha_-,\alpha_+],\\
2&\text{if}\ \, \theta>\theta_{\rm c} \ \, \text{and} \ \,
\alpha\in\{\alpha_{-},\alpha_+\},\\
3&\text{if}\ \,
\theta>\theta_{\rm c} \ \, \text{and} \ \, \alpha\in
(\alpha_-,\alpha_+),
\end{cases}
\end{equation*}
where $\theta_{\rm c}$ is given in \eqref{eq:theta'c} and
$\alpha_\pm=\alpha_\pm(\theta)$ are defined by~\eqref{eq:alpha_pm}.
\end{theorem}

Consider now the set of equations \eqref{uv-int}.
For each $m\in \{1,\dots,q-2\}$, consider
the functions
\begin{align}
\label{eq:L-intr} L_m(v;\theta):={}&(\theta-1)\mypp
\bigl(v^{k-1}+\dots+v\bigr)-m v^{k}
-(q-1-m),\\
\label{eq:K-intr} K_m(v;\theta):={}&
\frac{\bigl(v^{k-1}+\dots+v+1\bigr)^k
L_m(v;\theta)}{\bigl(v^{k-1}+\dots+v+L_m(v;\theta)\bigr)^k}.
\end{align}
For any
$\theta>1$ there is a unique value $v_m= v_m(\theta)>0$ such that
$$
L_m^*(\theta):=L_m(v_m;\theta)=\max_{v>0}L_m(v;\theta).
$$
Moreover, the function $\theta\mapsto L_m^*(\theta)$ is strictly
increasing. Denote by $\theta_m$ the unique value of $\theta>1$
such that
\begin{equation}\label{eq:Lmthetam}
L_m^*(\theta_m)=0.
\end{equation}
Thus, for any $\theta>\theta_m$ the range of the functions $v\mapsto
L_m(v;\theta)$ and $v\mapsto K_m(v;\theta)$ includes positive
values,
\begin{equation}\label{eq:V+}
\mathscr{V}_m^+(\theta):=\{v>0\colon L_m(v;\theta)>0\}\equiv
\{v>0\colon K_m(v;\theta)>0\}\ne \varnothing,
\end{equation}
and, therefore, the function
\begin{equation}\label{eq:alpham}
\alpha_m(\theta):=\frac{1}{\ln\theta} \max_{v\in
\mathscr{V}_m^+(\theta)}\ln K_m(v;\theta)=\frac{\ln
K_m^*(\theta)}{\ln\theta},\qquad \theta>\theta_m,
\end{equation}
is well defined, where
$$
K_m^*(\theta):=\max_{v\in \mathscr{V}_m^+(\theta)} K_m(v;\theta).
$$


\begin{theorem}\label{th:3.8} For each $m\in \{1,\dots,q-2\}$,
let\/ $\nu_{m}(\theta,\alpha)$ denote the number of TISGMs corresponding to positive
solutions $(u,v)$ of the system~\eqref{uv-int}. Then
$\nu_{m}(\theta,\alpha)\ge 1$ if and only if\/ $\theta>\theta_m$ and
$\alpha\le \alpha_m(\theta)$.
\end{theorem}

\strut{}In the case $q=3$, there is an additional
critical value
\begin{equation}\label{eq:theta-cr-q=3}
\tilde{\theta}_1=\tilde{\theta}_1(k):=\frac{5-k+\sqrt{49k^2+62k+49}}{6\myp(k-1)}.
\end{equation}

For $q\ge2$, consider the following subsets of the half-plane
$$\{\theta\ge1\}=\{(\theta,\alpha)\colon
\theta\ge1,\,\alpha\in\mathbb{R}\}.$$
\begin{align}
\notag A_q:={}&\{\theta>\theta_{\rm c},\,\alpha_{-}(\theta)\le
\alpha\le \alpha_{+}(\theta)\},\\[.2pc]
\label{eq:Bq} B_q:={}& \left\{
\begin{alignedat}{3}
& \,\varnothing     &&\text{if\ \,}q=2,\\
&\{\theta>\theta_1,\,\alpha<\alpha_1(\theta)\}
\cup\{\theta>\tilde{\theta}_1,
\,\alpha=\alpha_1(\theta)\}\quad&&\text{if\ \,}q=3,\\
&\{\theta>\theta_1,\,\alpha\le \alpha_1(\theta)\}\quad&&\text{if\
\,}q\ge4.
\end{alignedat}
\right.
\end{align}
Denote the total number of TISGMs corresponding to positive solutions
$\boldsymbol{z}=(z_1,\dots,z_{q-1})$ of the system of
equations~\eqref{pt1} by $\nu(\theta,\alpha)$ where $\theta\ge1$,
$\alpha\in\mathbb{R}$ (of course, this number also depends on $k$
and~$q$). Theorems \ref{th:3.7} and \ref{th:3.8} is
summarized as follows.

\begin{theorem}[Non-uniqueness]\label{th:non-U}\mbox{}
\begin{itemize}
\item[\rm(a)]
If\/ $q=2$ then $\nu(\theta,\alpha)\ge2$ if and only if
$(\theta,\alpha)\in A_2\cup B_2=A_2$.

\vspace{.3pc}
\item[\rm (b)] If\/ $q=3$ then
$\nu(\theta,\alpha)\ge 2$ if $(\theta,\alpha)\in A_3\cup B_3$. The\
``only if''\ statement holds true at least for $k=2,3,4$.

\vspace{.3pc}
\item[\rm (c)] If\/ $q\ge 4$ then
$\nu(\theta,\alpha)\ge 2$ if and only if $(\theta,\alpha)\in A_q\cup
B_q$.
\end{itemize}
\end{theorem}

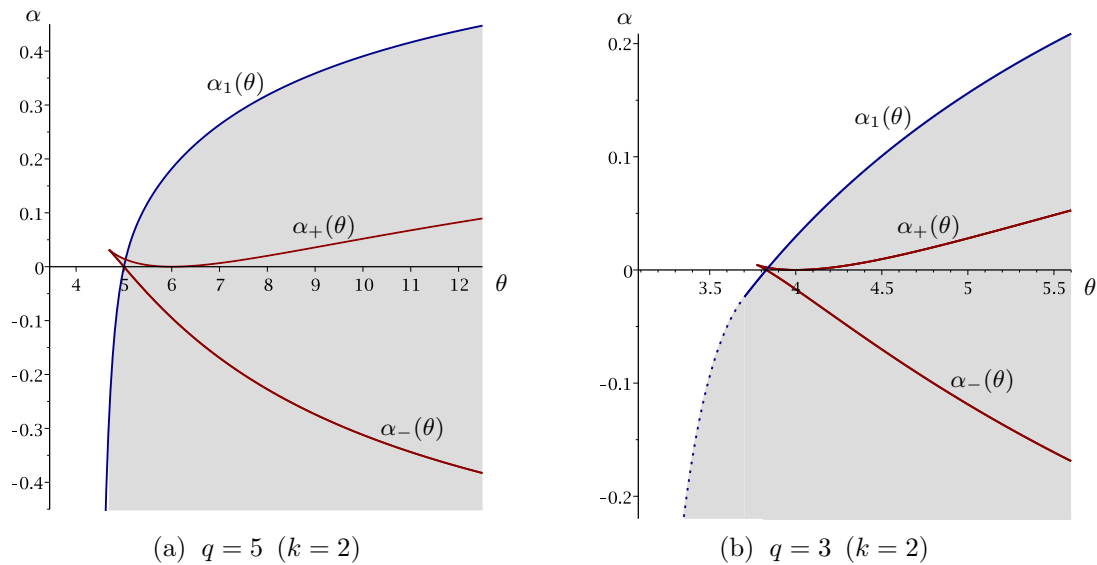
\begin{figure} \centering
\subfigure[\,$q=5$ \,($k=2$)]
{\hspace{-1pc}\includegraphics[width=6.7cm]{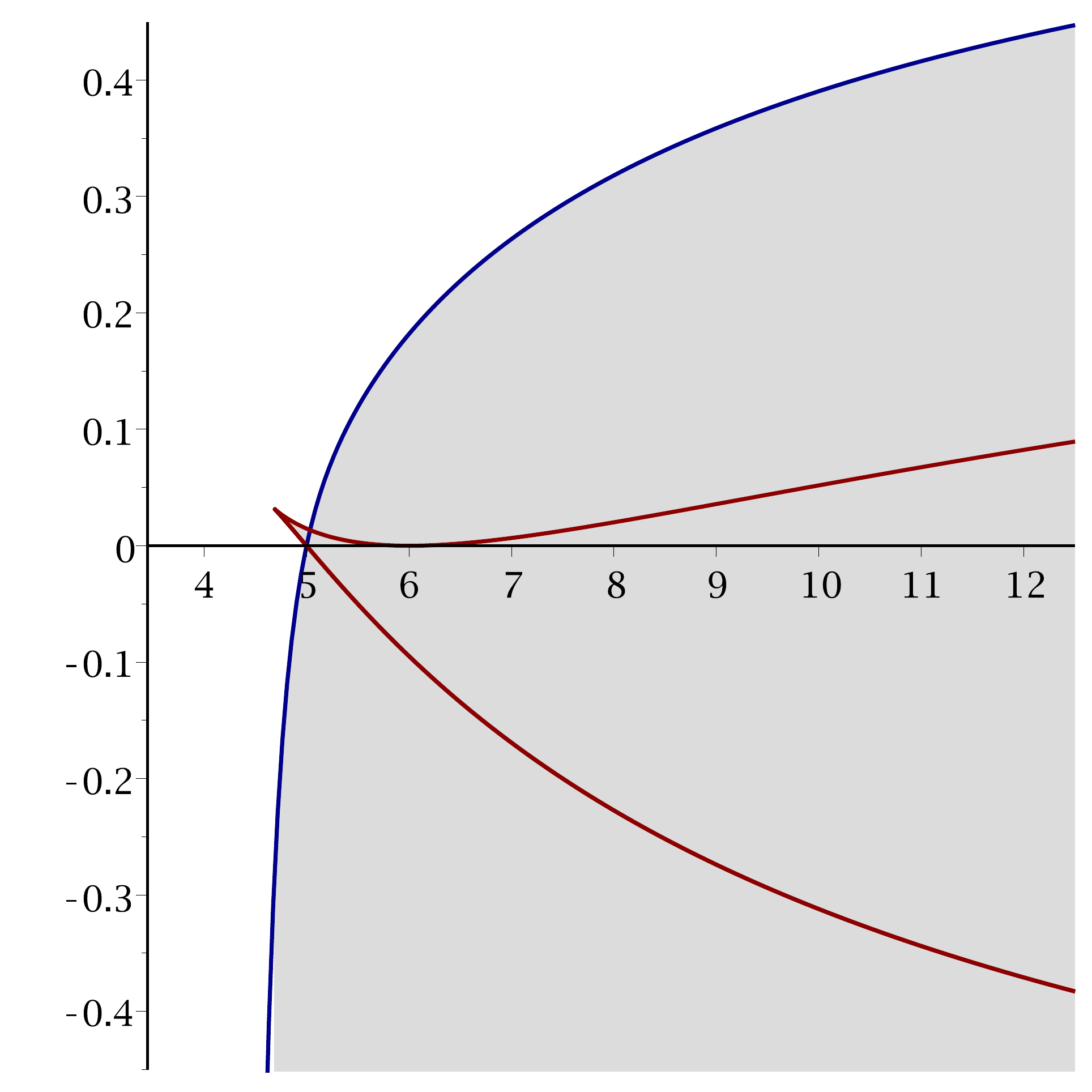}
\put(2,86){\mbox{\scriptsize$\theta$}}
\put(-174,188){\mbox{\scriptsize$\alpha$}}
\put(-106,162){\mbox{\scriptsize$\alpha_1(\theta)$}}
\put(-41,32){\mbox{\scriptsize$\alpha_{-}(\theta)$}}
\put(-74,108){\mbox{\scriptsize$\alpha_{+}(\theta)$}}
\label{Fig3a} } \hspace{1.5pc} \subfigure[\,$q=3$ \,($k=2$)]
{\raisebox{-.3pc}{\includegraphics[width=6.7cm]{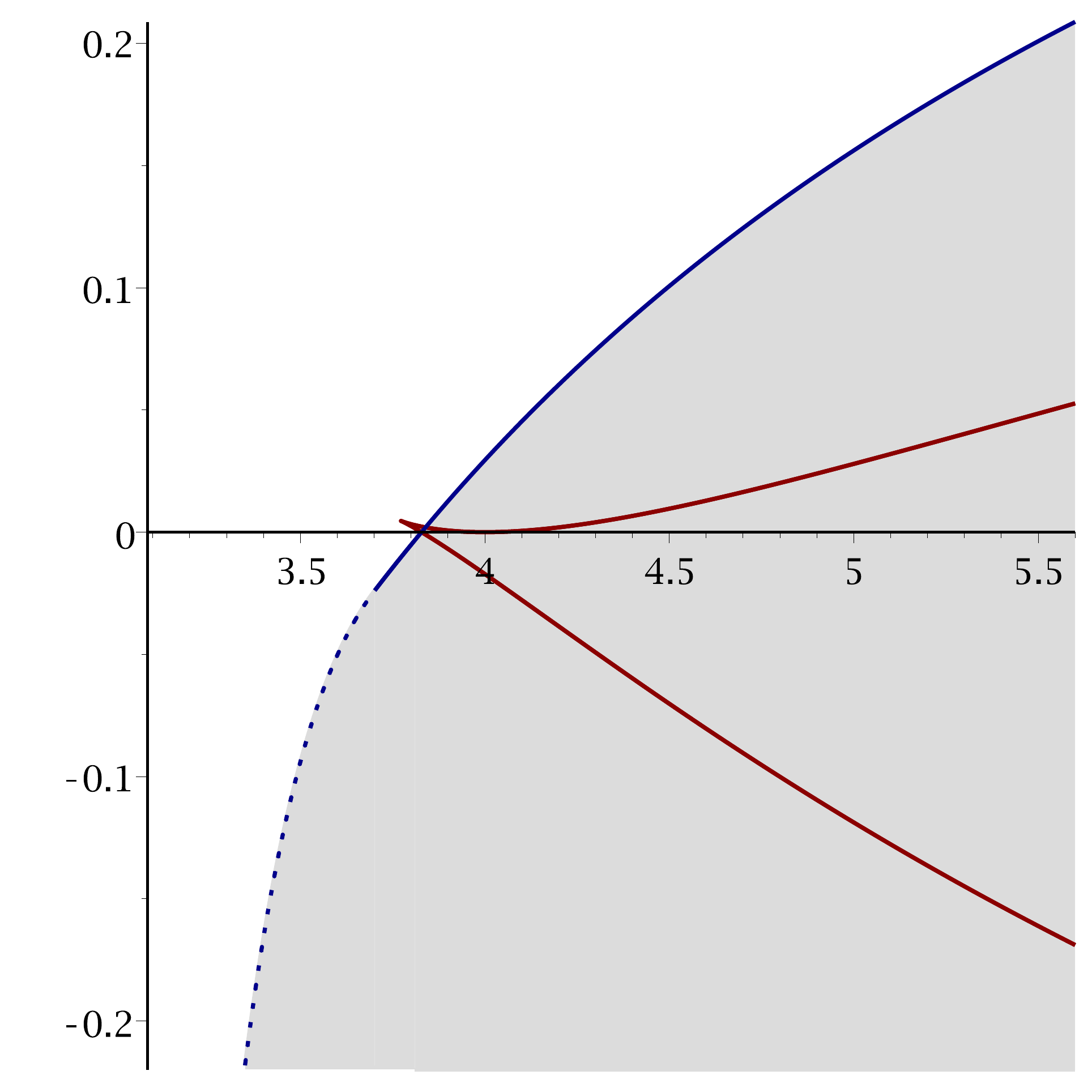}}\label{Fig3b}
\put(2,84.8){\mbox{\scriptsize$\theta$}}
\put(-173,187){\mbox{\scriptsize$\alpha$}}
\put(-84,149){\mbox{\scriptsize$\alpha_1(\theta)$}}
\put(-48,50){\mbox{\scriptsize$\alpha_{-}(\theta)$}}
\put(-68,108){\mbox{\scriptsize$\alpha_{+}(\theta)$}}}

\caption{The phase diagram for the Potts model~\eqref{eq:xi0}
showing the non-uniqueness region (shaded in grey) according to
Theorem~\ref{th:non-U}: (a) regular case $q\ge 4$ (shown here for
$q=5$); (b) special case $q=3$, both with $k=2$. The critical
boundaries are determined by (parts of) the graphs of the functions
$\alpha_\pm(\theta)$ and $\alpha_1(\theta)$ defined in
\eqref{eq:alpha_pm} and \eqref{eq:alpham}, respectively. The dotted
part of the boundary on panel (b), given by
$\alpha=\alpha_1(\theta)$, $\theta\in(\theta_1, \tilde{\theta}_1]$
(see \strut{}formula~\eqref{eq:Bq} with $q=3$), is excluded from the
shaded region (see Theorem~\ref{th:non-U}(b); here,
$\theta_1=3$ and
$\tilde{\theta}_1=\frac12\bigl(1+\sqrt{41}\myp\bigr)\doteq 3.7016$.
}\label{Fig3}
\end{figure}

Theorem \ref{th:non-U} provides a sufficient and (almost) necessary
condition for the uniqueness of solution of~\eqref{pt1}, illustrated
in Figure~\ref{Fig3} for $q=5$ and $q=3$, both with $k=2$.

\smallskip
To conclude this subsection, the following result describes a few
cases where it is possible to estimate the maximal number of
solutions of the system~\eqref{pt1}.
\begin{theorem}\label{th:k=2}\mbox{}
\begin{itemize}
\item[\rm (a)]
If $q=2$ then $\nu(\theta,\alpha)\le3$\textup{;} moreover,
$\nu(\theta,\alpha)=3$ for all $\theta\ge1$ large enough.

\smallskip
\item[\rm (b)] Let $\alpha=0$ and $k\ge2$. Then $\nu(\theta,0)\le 2^q-1$ for all
$\theta\ge1$\textup{;} moreover, $\nu(\theta,0)=2^q-1$ for all
$\theta$ large enough.

\smallskip
\item[\rm (c)] If $k=2$ then $\nu(\theta,\alpha)\le 2^q-1$ for all
$\theta\ge1$ and $\alpha\in\mathbb{R}$.
\end{itemize}
\end{theorem}

The general case $q\ge3$ with $\alpha\in\mathbb{R}$ was first
addressed by Peruggi et al.~\cite{Peruggi1} (and continued
in~\cite{Peruggi2}) using physical argumentation. In particular,
they correctly identified the critical point $\theta_{\rm c}$
\cite[equation~(22), page~160]{Peruggi2} (cf.~\eqref{eq:theta'c})
and also suggested an explicit critical boundary in the phase
diagram for $\alpha\ge0$, defined by the expression
\cite[equation~(21), page~160]{Peruggi2}:
$$
\tilde{\alpha}_{-}(\theta)=\frac{(k+1)\ln\bigl(1+(q-2)/\theta\bigr)-(k-1)\ln(q-1)}{2\ln\theta}.
$$

\section{Periodic SGMs}

\subsection{Definitions}

Let $\Gamma^k=(V, L)$ be a Cayley tree of order $k\geq 2$.

Consider $q$-state Potts model with an external field:
\begin{equation}\label{rr2}
H(\sigma)=-J\sum_{\langle x,y\rangle\in L}\delta_{\sigma(x)\sigma(y)}-\alpha\sum_{x\in V}
\delta_{1\sigma(x)}
\end{equation}
where $J, \alpha\in \mathbb R$.

By Theorem \ref{ep} we have that to each splitting Gibbs measures of this
model corresponds a set of vectors $h_x\in \mathbb R^{q-1}$, $x\in V$
which satisfies
\begin{equation}\label{rr5}
h_x=\sum_{y\in S(x)}F(h_y,\theta,\alpha),
\end{equation}
where $F: h=(h_1, \dots,h_{q-1})\in \mathbb R^{q-1}\to
F(h,\theta,\alpha)=(F_1,\dots,F_{q-1})\in \mathbb R^{q-1}$ defined as:
$$F_i=\alpha\beta\delta_{1i}+\ln\left({(\theta-1)e^{h_i}+\sum_{j=1}^{q-1}e^{h_j}+1\over
\theta+ \sum_{j=1}^{q-1}e^{h_j}}\right),$$
$\theta=\exp(J\beta)$, and $S(x)$ is the set of all direct
successors of $x\in V$ (see (\ref{sx})).

Recall $G_k$ is the group which represents the Cayley tree (see Section 2).
Let $G_k/G_k^*=\{H_1,...,H_r\}$ be the quotient group, where $G_k^*$ is a
normal subgroup of index $r\geq 1$.

\begin{definition}\label{dph} A set of vectors  $h=\{h_x,\, x\in G_k\}$
is said to be $G_k^*$-periodic, if \index{$G_k^*$-periodic}  $h_{yx}=h_x$ for any $ x\in
G_k$ and $y\in G^*_k.$ A $G_k-$ periodic set is called translation-invariant. \index{translation-invariant}
\end{definition}

\begin{definition}\label{dpw}  A set of vectors  $h=\{h_x,\, x\in G_k\}$
is said to be  $G_k^*$-weakly periodic, \index{$G_k^*$-weakly periodic} if $h_{x}=h_{ij}$ for
$x\in H_i$ and $x_\downarrow\in H_j$ for any $x\in G_k$.
\end{definition}

\begin{definition}\label{dpg} A measure $\mu$ is said to be
$G_k^*$-periodic (weakly periodic), if it corresponds to the
$G_k^*$-periodic (weakly periodic) set of vectors $h$. The $G_k-$
periodic measure is said to be translation-invariant.
\end{definition}

In this section we review the  results of  \cite{Hakimov}-\cite{KK}, \cite{RK}, \cite{RM}-\cite{Ras} related to (weakly) periodic
Gibbs measures of the Potts model.

Let $G^{(2)}_k$ be the subgroup in $G_k$ consisting of all words
of even length. Clearly, $G^{(2)}_k$ is a subgroup of index 2.

The following theorem characterizes all periodic Gibbs measures.

\begin{theorem}\label{rrt6} Let $K$ be a normal subgroup of finite
index in $G_k$. Then each $K$-periodic Gibbs measure for the Potts model (\ref{rr2})
 is either translation-invariant or $G^{(2)}_2$-periodic.
 \end{theorem}

 Since TISGMs are given in previous sections,
by Theorem \ref{rrt6} it suffices consider $G_k^{(2)}$-periodic
Gibbs measure only. Such measures correspond to the set of
vectors $h=\{h_x\in \mathbb R^{q-1}: \, x\in G_k\}$ of the form
$$h_x=\left\{%
\begin{array}{ll}
   {\bf h}^1, \ \ \ $ if $ |x| \ \ \mbox{is even} $,$ \\
    {\bf h}^2, \ \ \ $ if $ |x| \ \ \mbox{is odd} $,$
\end{array}%
\right. $$
where $|x|$ denotes the length of the word $x$ and the vectors
${\bf h}^1=(h_1,h_2,...,h_{q-1})$,  ${\bf h}^2=(l_1,l_2,...,l_{q-1})$,  by (\ref{rr5}), should
satisfy
$$
\left\{
\begin{array}{ll}
    h_{i}=k\ln\cfrac{(\theta-1)\exp(l_i) + \sum_{j=1}^{q-1}\exp({l_j})+1}{\sum_{j=1}^{q-1}\exp({l_j})+\theta},\\[5 mm]
    l_{i}=k\ln\cfrac{(\theta-1)\exp(h_i) + \sum_{j=1}^{q-1}\exp({h_j})+1}{\sum_{j=1}^{q-1}\exp({h_j})+\theta}
\end{array}
i=1, \dots, q-1. \right.$$\

Denoting $\exp(h_i)=x_i,\ \exp(l_i)=y_i$ the last system can be written as
\begin{equation}\label{kr4}
\begin{split}
    x_{i}=\left({(\theta-1)y_i + \sum_{j=1}^{q-1}y_j+1\over \sum_{j=1}^{q-1}y_j+\theta}\right)^k, \\[4 mm]
    y_{i}=\left({(\theta-1)x_i + \sum_{j=1}^{q-1}x_j+1\over \sum_{j=1}^{q-1}x_j+\theta}\right)^k
    \end{split}%
    i=1, \dots, q-1.
\end{equation}

Consider mapping $W:\mathbb R_+^{q-1}\times \mathbb R_+^{q-1} \rightarrow
\mathbb R_+^{q-1}\times \mathbb R_+^{q-1},$ defined by
\begin{equation}\label{kr5}
\begin{split}
    x_{i}'=\left({(\theta-1)y_i + \sum_{j=1}^{q-1}y_j+1\over \sum_{j=1}^{q-1}y_j+\theta}\right)^k, \\[4 mm]
    y_{i}'=\left({(\theta-1)x_i + \sum_{j=1}^{q-1}x_j+1\over \sum_{j=1}^{q-1}x_j+\theta}\right)^k
    \end{split}%
    i=1, \dots, q-1.
\end{equation}

Note that (\ref{kr4}) is the equation $z=W(z)$, with
$$z=(x_1, \dots, x_{q-1}, \, y_1, \dots, y_{q-1})$$
that is equation for fixed points of the mapping $W$.

For $m\leq q-1$ introduce the set
$$I_m=\{x=(x_1,\dots, x_{q-1})\in \mathbb R_+^{q-1}: x_1=\dots=x_m, \ \ x_{m+1}=\dots=x_{q-1}=1\}.$$
Let  $S_{q-1}$ be the group of permutations on $\{1,\dots, q-1\}$.

For $\pi=(\pi(1), \dots, \pi(q-1))\in S_{q-1}$ and $x\in \mathbb R^{q-1}$
define
$$\pi x=(x_{\pi(1)}, \dots , x_{\pi(q-1)}),$$
and for a subset $I\subset \mathbb R^{q-1}$ define
$$\pi I=\{\pi x: x\in I\}.$$

The following lemma is useful to find special solutions of (\ref{kr4}).
\begin{lemma}\label{he} For any $\pi\in S_{q-1}$ and $m\leq q-1$  the set
$\pi I_m\times \pi I_m$ is invariant with respect to the mapping $W$.
\end{lemma}

 \subsection{Zero external field}\
In this section consider the case $\alpha=0$.\\

{\it Ferromagnetic case.} \ \ In \cite{RMM},
it is shown that in the ferromagnetic case (i.e. $J>0$) each periodic measure is translation-invariant:

\begin{theorem} Let $k\geq2, \ q\geq 2, \ J>0, \ \alpha=0$.
Then any periodic Gibbs measure of the model (\ref{rr2}) is translation-invariant.
\end{theorem}

 For $q=3$ this theorem was proved in \cite{RK}.\\

{\it Anti-ferromagnetic case.} \ \
Note that the case $q=2$
corresponds to anti-ferromagnetic Ising model, which has periodic measures for any $k\geq 2$ (\cite{Ge}, \cite{R}).

In this subsection we consider the anti-ferromagnetic case (i.e. $J<0$) and review results of   \cite{Hakimov1}, \cite{KK} and \cite{RMM}   that for
$k=2$, $q\geq 3$ there is no periodic (except translation-invariant) Gibbs measures.  But for
$k\geq 3$, $q\geq 2$ there are several periodic (non-translation-invariant) Gibbs measures.


\begin{theorem} Let $k=2$, $q\geq 3$, $J<0$, and $\alpha=0$.
Then the $G_k^{(2)}$-periodic Gibbs measure for the Potts
model is unique. Moreover, this measure coincides with the unique translation-invariant Gibbs measure.
\end{theorem}


In case $k\geq 3$, we have
the following

\begin{theorem}\label{davr} Let  $k\geq 3$, $q\geq 3$, $J<0$ and $\theta_{cr}={k-q+1\over k+1}$.
 Then for the Potts model for $0<\theta<\theta_{cr}$ there are exactly three $G_k^{(2)}-$
periodic Gibbs measures, corresponding to solutions on $I_m\times I_m$.
Moreover, only one of these measures is translation-invariant.
\end{theorem}

 Using symmetry of solutions, one obtains the following:

\begin{theorem} For the Potts model on the Cayley tree of order  $k\geq 3$, if $3\leq
q<k+1$ and  $0<\theta<\theta_{cr}$ then the number of $G_k^{(2)}-$periodic Gibbs measures is at least
$2\cdot (2^{q}-1)$.
\end{theorem}

\subsection {Non-zero external field.}

In case $\alpha\ne 0$, the system of equations for periodic solutions becomes
more complicated to solve. Because in this case the symmetry of solutions will be lost.

In this subsection following \cite{RK} we consider a particular case of parameters and
give a result that under some conditions on parameters the model has periodic (non-translation-invariant) Gibbs measures.

Let  $k=2$, $q=3$ and  $e^{\alpha \beta}=\lambda,$ where $\alpha\in \mathbb R $
is the external field, and as before, $\beta=\frac{1}{T},$ $T>0$.

Periodic measures correspond to positive solutions of
%
\begin{equation}\label{kr19}
\left\{%
\begin{array}{ll}
    x=\zeta \cdot \cfrac{\theta y^2+2}{y^2 +\theta+1}\\[2mm]
    y=\zeta\cdot \cfrac{\theta x^2+2}{x^2 +\theta+1}.
    \end{array}%
\right.
\end{equation}
where $\zeta=\sqrt{\lambda}.$

 This system, in case $x=y$, can be written as
\begin{equation}\label{kr20}x^3-\zeta \theta x^2+(\theta+1)x-2\zeta=0,
\end{equation}
which gives translation-invariant measures.

Substituting the expression for $y$ from the second
equation in the first equation, we obtain a fifth-degree equation in $x$.
Dividing the fifth-degree polynomial
from the obtained equation by the third-degree polynomial from (\ref{kr20}), we obtain the quadratic equation
\begin{equation}\label{kr21}(\zeta^2\theta^2+\theta+1)x^2+(-2\zeta+\zeta\theta(\theta+1))x+2\zeta^2\theta+(\theta+1)^2=0,
\end{equation}
which describes periodic Gibbs measures. We now study solutions of this quadratic equation. For this, we
calculate its discriminant taking into account that $\sqrt{\lambda}=\zeta$:
$$D=-8\theta^3\lambda^2-[3\theta^2(\theta+1)^2+12\theta(\theta+1)-4]\lambda-4(\theta+1)^3.$$
It is known that equation (\ref{kr21}) has two positive real roots if $D>0$ and if
$-2\zeta+\zeta\theta(\theta+1)=\zeta(\theta^2+\theta-2)<0.$
 Solving these inequalities, we find that equation (\ref{kr21})
 has two real solutions if $0<\theta< \frac
{\sqrt{17}-3}{6}=\theta^{*}$ and
$\lambda_1(\theta)<\lambda<\lambda_2(\theta),$ where
$\lambda_1(\theta), \ \lambda_2(\theta)$ has the following form
\begin{equation}\label{kr23}
\begin{split}
\lambda_1(\theta)=\frac{-3[3\theta^2(\theta+1)^2+12\theta(\theta+1)-4]-
(\theta-1)(\theta+2)\sqrt{K}}{48\theta^3},\\[3mm]
\lambda_2(\theta)=\frac{-3[3\theta^2(\theta+1)^2+12\theta(\theta+1)-4]
+(\theta-1)(\theta+2)\sqrt{K}}{48\theta^3},
\end{split}
\end{equation}
where
$$K=(\theta-1)(\theta+2)(9\theta^2+9\theta-2).$$
Note that for $0<\theta<\theta^{*}$ both
$\lambda_1(\theta)$ and $\lambda_2(\theta)$ are positive.\

If $\theta^{*}<\theta<1,$ then equation (\ref{kr21}) has no real solutions.
Moreover, equation (\ref{kr21}) also has no real solutions for $\theta>1.$

Let now $\theta=\theta^{*}.$ Then for $\lambda=\lambda^{*}$, where
\begin{equation}\label{kr24}
\lambda^{*}=-\frac {3(\theta^{*})^2(\theta^{*}+1)^2+12\theta^{*}(\theta^{*}+1)-4}{16(\theta^{*})^3},
\end{equation}
the discriminant of the equation (\ref{kr21}) is $D=0.$ Consequently, the equation (\ref{kr21}) has unique solution:
$$x^{*}=\frac{\zeta^{*}(2-\theta^{*}-(\theta^{*})^2)}{2((\zeta^{*})^2(\theta^{*})^2+\theta^{*}+1)},$$
where $\zeta^{*}=\sqrt{\lambda^{*}}.$


Thus we have the following theorem:

\begin{theorem} For the Potts model, in case  $k=2$, $q=3$, $\alpha\ne 0$ we have
\begin{itemize}
\item[1.] \ \ If $\lambda\in(\lambda_1(\theta), \
\lambda_2(\theta))$ and $\theta\in(0, \ \theta^{*}),$ then there are at least two $G_k^{(2)}-$ periodic
(non-translation-invariant) Gibbs measures,
 where
$\theta^{*}=\frac {\sqrt{17}-3}{6},$ and $\lambda_1, \ \lambda_2$
defined by (\ref{kr23});
\item[2.] \ \  For $\lambda=\lambda^{*}$ and  $\theta=\theta^{*},$ there are at least one
$G_k^{(2)}-$ periodic
(non-translation-invariant) Gibbs measure,
 where
$\lambda^{*}$ is defined in (\ref{kr24});
\end{itemize}
\end{theorem}

 \section{Weakly periodic measures.}

In this section following \cite{RM}-\cite{RMb} we give some weakly
periodic Gibbs measures (see Definition \ref{dpg}) of the Potts model.

The level of difficulty in describing of weakly periodic Gibbs
measures is related to the structure and index of the normal
subgroup relative to which the periodicity condition is imposed.

From Chapter 1 of \cite{R} (see also \cite{NR}) we know that in the group $G_k$, there is no normal
subgroup of odd index different from one. Therefore, we consider
normal subgroups of even indices. In this section we restrict ourself to the
case of indices two.

It is known (see  \cite{GRU} and \cite{R}) that any normal subgroup of index two  of the group
$G_k$ has the form
$$H_A=\left\{x\in G_k:\sum\limits_{i\in A}\omega_x(a_i){\rm -even} \right\},$$
where $\emptyset \neq A\subseteq N_k=\{1,2,\dots,k+1\}$, and
$\omega_x(a_i)$ is the number of letters $a_i$ in a word $x\in
G_k$.

Let $A\subseteq N_k$ and $H_A$ be the corresponding normal
subgroup of index two. We note that in the case $|A|=k+1$, i.e.,
in the case $A = N_k$, weak periodicity coincides with ordinary
periodicity.

Consider $G_k/H_A=\{H_A,G_k\setminus H_A\}$ the quotient group.

For simplicity of notations, we denote $H_0=H_A,$ $ H_1=G_k\setminus H_A$.

The $H_0$ -weakly periodic collections of vectors $h=\{h_x\in
\mathbb{R}^{q-1}: \, x\in G_k\}$ have the form:

\begin{equation}\label{wwp}
h_x=\left\{%
\begin{array}{ll}
    h_1, & \textrm{if} \ \ {x_{\downarrow} \in H_0}, \ x \in H_0 \\
    h_2, & \textrm{if} \ \ {x_{\downarrow} \in H_0}, \ x \in H_1 \\
    h_3, & \textrm{if} \ \ {x_{\downarrow} \in H_1}, \ x \in H_0 \\
    h_4, & \textrm{if} \ \ {x_{\downarrow} \in H_1}, \ x \in H_1.
\end{array}%
\right.
\end{equation}
Here $h_i=(h_{i1},h_{i2},...,h_{iq-1}),$ $i=1,2,3,4.$
\begin{remark} The definition of $x_{\downarrow}$ depends on the (fixed) root $e$ of the Cayley tree.
Therefore the function $h_x$ given in (\ref{wwp}) is defined for $x\ne e$.
\end{remark}

By (\ref{rr5})
we then have
\begin{equation}\label{mn5}
\left\{%
\begin{array}{ll}
    h_{1}=(k-|A|)F(h_{1},\theta,\alpha)+|A|F(h_{2},\theta,\alpha) \\[2mm]
    h_{2}=(|A|-1)F(h_{3},\theta, \alpha)+(k+1-|A|)F(h_{4},\theta, \alpha) \\[2mm]
    h_{3}=(|A|-1)F(h_{2},\theta, \alpha)+(k+1-|A|)F(h_{1},\theta, \alpha) \\[2mm]
    h_{4}=(k-|A|)F(h_{4},\theta, \alpha)+|A|F(h_{3},\theta, \alpha).
\end{array}%
\right.\end{equation}

Let us give an  examples geometrically presenting
the vector-values of $h_x$ on the Cayley tree.

\begin{example} Let $A=\{1\}$. Then
$$H_A=H_{\{1\}}=\left\{x\in G_k: \ \ {\rm number \, of} \ \ a_1 \ \ {\rm in} \ \ x \ \ {\rm is \, even} \right\},$$
In this case $H_0=H_{\{1\}}$ (red points in Fig. \ref{fwp}) and $H_1=G_k\setminus H_0$
is the set of words with an odd number of $a_1$ (green points in the Fig. \ref{fwp}).
\begin{figure}
\vspace{-.5pc} \centering
  \includegraphics[width=12cm]{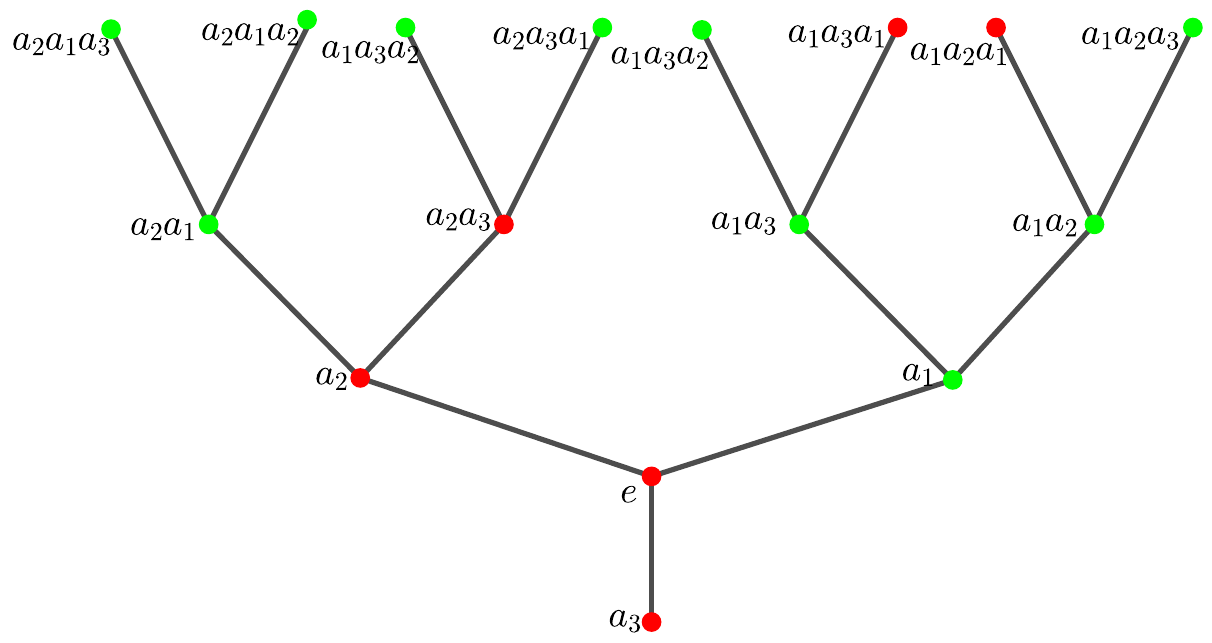}
  \caption{Partition (red and green) of the Cayley tree with respect to
  subgroup $H_{\{1\}}$ of index two.}\label{fwp}
\end{figure}

In the Fig.\ref{fww} we showed the distribution of the four values of function (\ref{wwp}).
\begin{figure}
\vspace{-.5pc} \centering
  \includegraphics[width=12cm]{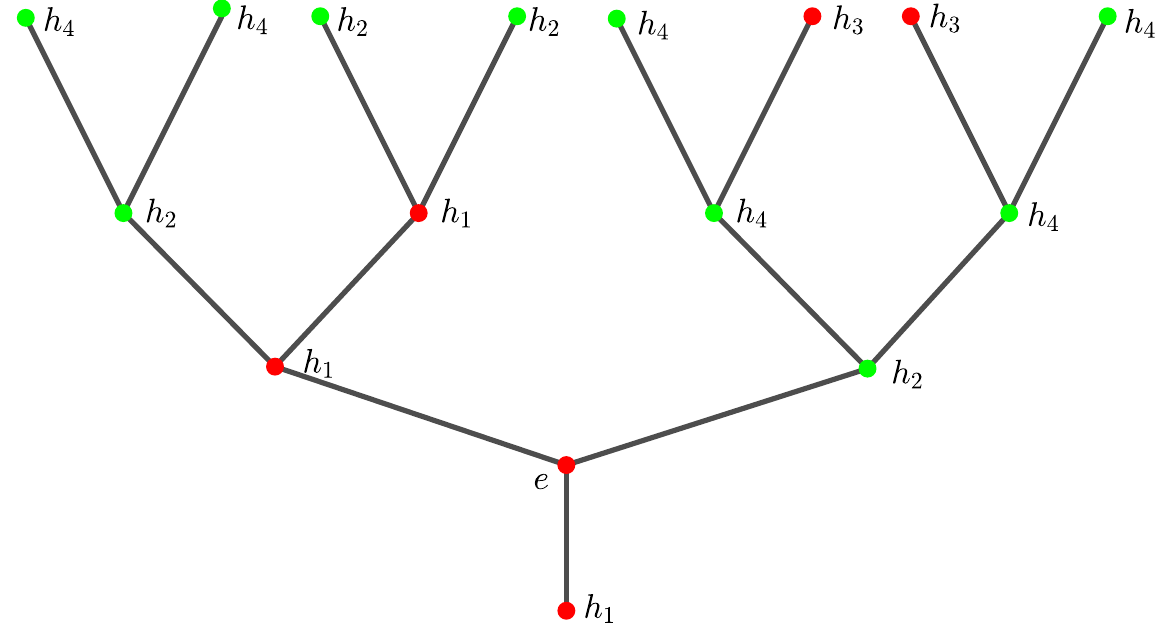}
  \caption{The values of the weakly-periodic function (\ref{wwp}) in case $A=\{1\}$, i.e. subgroup $H_{\{1\}}$ of index two.}\label{fww}
\end{figure}
\end{example}

Denote
$z_{ij}=\exp{h_{ij}},\lambda_j=\exp(\alpha_j), i=1,2,3,4,
j=1,2,...q-1$. Then the last system of equations can be rewritten
as
\begin{equation}\label{mn6}
\left\{%
\begin{array}{ll}
    z_{1j}=\lambda_j \left({(\theta-1) z_{1j} + \sum_{i=1}^{q-1}z_{1i}+1\over \sum_{i=1}^{q-1}z_{1i}+\theta}\right)^{k-|A|}\left({(\theta-1) z_{2j} + \sum_{i=1}^{q-1}z_{2i}+1\over \sum_{i=1}^{q-1}z_{2i}+\theta}\right)^{|A|}\\[3mm]
    z_{2j}=\lambda_j \left({(\theta-1) z_{3j} + \sum_{i=1}^{q-1}z_{3i}+1\over \sum_{i=1}^{q-1}z_{3i}+\theta}\right)^{|A|-1}\left({(\theta-1) z_{4j} + \sum_{i=1}^{q-1}z_{4i}+1\over \sum_{i=1}^{q-1}z_{4i}+\theta}\right)^{k+1-|A|}\\[3mm]
    z_{3j}=\lambda_j \left({(\theta-1) z_{2j} + \sum_{i=1}^{q-1}z_{2i}+1\over \sum_{i=1}^{q-1}z_{2i}+\theta}\right)^{|A|-1} \left({(\theta-1) z_{1j} + \sum_{i=1}^{q-1}z_{1i}+1\over \sum_{i=1}^{q-1}z_{1i}+\theta}\right)^{k+1-|A|}\\[3mm]
    z_{4j}=\lambda_j \left({(\theta-1) z_{4j} + \sum_{i=1}^{q-1}z_{4i}+1\over \sum_{i=1}^{q-1}z_{4i}+\theta}\right)^{k-|A|} \left({(\theta-1) z_{3j} + \sum_{i=1}^{q-1}z_{3i}+1\over \sum_{i=1}^{q-1}z_{3i}+\theta}\right)^{|A|},
\end{array}%
\right.\end{equation} here $j=1,2,3,..,q-1.$

Define the map $K:\mathbb{R}^{4(q-1)}\rightarrow \mathbb{R}^{4(q-1)},$ as
\begin{equation}\label{mn6'}
\left\{%
\begin{array}{ll}
    z'_{1j}=\lambda_j \left({(\theta-1) z_{1j} + \sum_{i=1}^{q-1}z_{1i}+1\over \sum_{i=1}^{q-1}z_{1i}+\theta}\right)^{k-|A|}\left({(\theta-1) z_{2j} + \sum_{i=1}^{q-1}z_{2i}+1\over \sum_{i=1}^{q-1}z_{2i}+\theta}\right)^{|A|}\\[3mm]
    z'_{2j}=\lambda_j \left({(\theta-1) z_{3j} + \sum_{i=1}^{q-1}z_{3i}+1\over \sum_{i=1}^{q-1}z_{3i}+\theta}\right)^{|A|-1}\left({(\theta-1) z_{4j} + \sum_{i=1}^{q-1}z_{4i}+1\over \sum_{i=1}^{q-1}z_{4i}+\theta}\right)^{k+1-|A|}\\[3mm]
    z'_{3j}=\lambda_j \left({(\theta-1) z_{2j} + \sum_{i=1}^{q-1}z_{2i}+1\over \sum_{i=1}^{q-1}z_{2i}+\theta}\right)^{|A|-1} \left({(\theta-1) z_{1j} + \sum_{i=1}^{q-1}z_{1i}+1\over \sum_{i=1}^{q-1}z_{1i}+\theta}\right)^{k+1-|A|}\\[3mm]
    z'_{4j}=\lambda_j \left({(\theta-1) z_{4j} + \sum_{i=1}^{q-1}z_{4i}+1\over \sum_{i=1}^{q-1}z_{4i}+\theta}\right)^{k-|A|} \left({(\theta-1) z_{3j} + \sum_{i=1}^{q-1}z_{3i}+1\over \sum_{i=1}^{q-1}z_{3i}+\theta}\right)^{|A|},
\end{array}%
\right.\end{equation} here $j=1,2,3,..,q-1.$

Denote
\begin{equation}\label{mn7} I_m=\{(z_1,...,z_{q-1})\in
\mathbb{R}^{q-1}:z_1=...=z_m,
z_{m+1}=...=z_{q-1}=1\},\end{equation}
\begin{equation}\label{mn7'} M_m=\{(z^{(1)}, z^{(2)}, z^{(3)}, z^{(4)})\in
\mathbb{R}^{4(q-1)}:z^{(i)}\in I_m, i=1, 2, 3, 4\}, \end{equation}
here $m=1,2,...,q-1.$

 \begin{lemma}\label{lemma1} We have
 \begin{itemize}
 \item[1)] \ \ For any fixed $m\geq 1$, $\lambda>0$ and
\begin{equation}\label{lam}
\lambda_i=\left\{%
\begin{array}{ll}
\lambda \ \ \mbox{if} \ \ 1\leq i\leq m\\
1 \ \ \mbox{if} \ \ m<i\leq q-1
\end{array}%
\right.\end{equation} the set $M_m$ is an invariant set with respect to the map $K$, i.e. $K(M_m)\subset M_m$.
\item[2)] \ \  For $\alpha=0$ the sets $M_m$ are  invariant with respect to
the map $K$ for all $m=1,2,...,q-1.$
\end{itemize}
\end{lemma}

\textbf{Case $\alpha \neq 0$.} Let us consider the case
$\alpha\neq 0$ and $\lambda_i$ given by (\ref{lam}). For $\textbf{z}\in M_m$ we denote $z_i=z_{ij},
i=1,2,3,4; j=1,2,...,m.$. Then on the invariant set $M_m$ the system
of equations (\ref{mn6}) has the following form

\begin{equation}\label{mn8}
\left\{%
\begin{array}{ll}
    z_{1}=\lambda\left({(\theta +m-1) z_1+q-m\over {mz_1+\theta+q-m-1}}\right)^{k-|A|} \left({(\theta +m-1) z_2+q-m\over {mz_2+\theta+q-m-1}}\right)^{|A|}\\[2mm]
    z_{2}=\lambda\left({(\theta +m-1) z_3+q-m\over {mz_3+\theta+q-m-1}}\right)^{|A|-1} \left({(\theta +m-1) z_4+q-m\over {mz_4+\theta+q-m-1}}\right)^{k+1-|A|}\\[2mm]
    z_{3}=\lambda\left({(\theta +m-1) z_2+q-m\over {mz_2+\theta+q-m-1}}\right)^{|A|-1} \left({(\theta +m-1) z_1+q-m\over {mz_1+\theta+q-m-1}}\right)^{k+1-|A|}\\[2mm]
    z_{4}=\lambda\left({(\theta +m-1) z_4+q-m\over {mz_4+\theta+q-m-1}}\right)^{k-|A|} \left({(\theta +m-1) z_3+q-m\over {mz_3+\theta+q-m-1}}\right)^{|A|}.
\end{array}%
\right.\end{equation}

We introduce the notation
$$g_m(z)={(\theta +m-1) z+q-m\over {mz+\theta+q-m-1}}.$$

It is easy to prove the following

\begin{lemma}\label{lemma2} The function $g_m(z)$ is strictly decreasing
 for $0<\theta<1$, $1\leq m\leq q-1$ and strictly increasing for
 $1<\theta$.\end{lemma}

\begin{proposition}\label{pro1} Let $\textbf{z}=(z_1, z_2, z_3, z_4)$ is a
solution of the system of equations (\ref{mn8}). If $z_i=z_j$ for
some $i\neq j$, then $z_1=z_2=z_3=z_4.$
\end{proposition}

Consider now the anti-ferromagnetic Potts model (i.e. $0<\theta<1$).

Let $|A|=k$. Then system of equation
(\ref{mn8})
%
can be reduced to analyzing of the following system of
equations
\begin{equation}\label{mn13}
\left\{%
\begin{split}
    z_{2}=\lambda\left(g_m(z_3)\right)^{k-1}\cdot g_m(\lambda(g_m(z_3))^k)\\[2mm]
    z_{3}=\lambda\left(g_m(z_2)\right)^{k-1}\cdot g_m(\lambda(g_m(z_2))^k).
\end{split}%
\right.
\end{equation}
Introducing the notation
\begin{equation}\label{mn14}
\psi(z)=\lambda\left(g_m(z)\right)^{k-1}\cdot
g_m(\lambda(g_m(z))^k).
\end{equation}
Then we reduce the system of equations (\ref{mn13}) to the form
\begin{equation}\label{mn15}
\left\{%
\begin{array}{ll}
    z_{2}=\psi(z_3)\\
    z_{3}=\psi(z_2).
\end{array}%
\right.
\end{equation}
The number of solutions of the system (\ref{mn15}) coincides with the
number of solutions of the equation $\psi(\psi(z))=z$.

%

It is known (Theorem \ref{antiT}, see also \cite[page 109]{R})  that for anti-ferromagnetic case,
there exits a unique translation-invariant Gibbs measure corresponding to the
unique solution of the equation $z=\lambda g_m^k(z)$. We let $z_*$
denote this solution.

\begin{proposition}\label{pro2} For $k\geq 6$ and
$\lambda\in(\lambda_{c_1}, \lambda_{c_2})$ the system of equations
({\ref{mn15}}) has three solutions $(z_*, z_*), (z^*_2, z^*_3),
(z^*_3, z^*_2)$, where $\lambda_{c_i}= b^k_{i}, i=1,2$ and
\begin{equation}\label{mn17}
\begin{split}
b_1=\frac{(k-1-\sqrt{k^2-6k+1}){(1-\theta)(\theta+q-1)}z_*^{\frac{k-1}{k}}}{2(mz_*+\theta+q-m-1)^2},\\[2mm]
b_2=\frac{(k-1+\sqrt{k^2-6k+1}){(1-\theta)(\theta+q-1)}z_*^{\frac{k-1}{k}}}{2(mz_*+\theta+q-m-1)^2}.
\end{split}
\end{equation}\end{proposition}

We have thus the following theorem.

\begin{theorem}\label{thm3} Let $|A|=k, k\geq 6$ and
$\lambda\in(\lambda_{c_1}, \lambda_{c_2})$. Then for the
anti-ferromagnetic Potts model with the special external field
(given by (\ref{lam})) there are at last two $H_A-$ weakly
periodic (non-periodic) Gibbs measures, where $\lambda_{c_i}=
b^k_{i}, i=1,2$.
\end{theorem}

\textbf{Case $\alpha=0$.} In this case the system of equations
(\ref{mn6}) on the invariant set $M_m, m=1,2,...,q-1$ can be reduced to the
following system of equations:

\begin{equation}\label{mn16'}
\left\{%
\begin{array}{ll}
    z_{1}=\left({(\theta+m-1) z_1+q-m\over {mz_1+\theta+q-m-1}}\right)^{k-|A|}\left({(\theta+m-1) z_2+q-m\over {mz_2+\theta+q-m-1}}\right)^{|A|}\\[2mm]
    z_{2}=\left({(\theta+m-1) z_3+q-m\over {mz_3+\theta+q-m-1}}\right)^{|A|-1} \left({(\theta+m-1) z_4+q-m\over {mz_4+\theta+q-m-1}}\right)^{k+1-|A|}\\[2mm]
    z_{3}=\left({(\theta+m-1) z_2+q-m\over {mz_2+\theta+q-m-1}}\right)^{|A|-1} \left({(\theta+m-1) z_1+q-m\over {mz_1+\theta+q-m-1}}\right)^{k+1-|A|}\\[2mm]
    z_{4}=\left({(\theta+m-1) z_4+q-m\over {mz_4+\theta+q-m-1}}\right)^{k-|A|}\left({(\theta+m-1) z_3+q-m\over {mz_3+\theta+q-m-1}}\right)^{|A|}.
\end{array}%
\right.\end{equation}

The following proposition is similar to Proposition \ref{pro1}.

\begin{proposition}\label{pro3} Let $m\in\{1,2,...,q-1\}$
be fixed and $\textbf{z}=(z_1, z_2, z_3, z_4)$ is a solution of
the system of equations (\ref{mn16'}). If $z_i=z_j$ for some $i\neq
j$, then $z_1=z_2=z_3=z_4.$
\end{proposition}

\begin{theorem}\label{thm4} Let $|A|=k$ and $k\geq 6$. If
one of the following conditions is satisfied

1) $\frac{4k}{k+1+\sqrt{k^2-6k+1}}\leq q
<\frac{4k}{k+1-\sqrt{k^2-6k+1}}$ and $0<\t<\t_2$;

2) $ q \leq \frac{4k}{k+1+\sqrt{k^2-6k+1}}$ and $\t_1<\t<\t_2,$

then there are at last $2^q-2$ weakly periodic (non-periodic)
Gibbs measures, where
\begin{equation}\label{mn*15}
\begin{split}
\t_1={\frac {4\,k-kq-q- q\sqrt{{k}^{2}-6\,k+1}}{4k}},\\[2mm]
\t_2={\frac{4\,k-kq-q+q\sqrt{{k}^{2}-6\,k+1}}{4k}}.
\end{split}
\end{equation}
\end{theorem}

\begin{remark}\label{rk1}
The new Gibbs measures described in Theorem \ref{thm3} and Theorem
\ref{thm4} allow to describing a continuous set of non-periodic
Gibbs measures (see Chapter 2 of \cite{R} and \cite{Ru} for such  constructions) different from the previously known ones.\end{remark}

\begin{remark}\label{rk1} For the case $m=q-1$ (see (\ref{mn7'})),  Theorem \ref{thm4} coincides with Theorem 3 in
\cite{RM}.\end{remark}

{\bf References and comments.}
This review is based on works \cite{BR}, \cite{FV}, \cite{GRRR}, \cite{GRR}, \cite{GRU},
\cite{Ga8}, \cite{Ga9}, \cite{Ge}, \cite{HaKu04}, \cite{YH}, \cite{Hakimov}-\cite{KK}, \cite{Ke},
\cite{KRK}, \cite{KR}, \cite{MarJ}, \cite{MSW}, \cite{Mos}, \cite{RM}-\cite{Ras}, \cite{Ro9}, \cite{R}-\cite{Rc1},
\cite{RR}, \cite{RKK}. Many of these papers written after 2013 and are related to the Potts model on trees.

In my opinion, one of nice results
of the theory of Gibbs measures on trees is Theorem \ref{Theorem4},
 which gives full description of TISGMs for $q$-state Potts model.
 Namely, it states that there are $[q/2]$ distinct critical temperatures changing the number
 of TISGMs and there can be up to $2^q-1$ such measures.
Before of this result it was only known one critical temperature and $q$ TISGMs \cite{Ga8},\cite{Ga9}.

Applications of results of \cite{MSW} allowed to give conditions of extremality of
all TISGMs. Remaining results of this review are related to recently obtained results on
boundary configurations, non-zero external fields,
periodic and weakly periodic Gibbs measures and free energies of the Potts model.

For future reading see \cite{Aq}, \cite{BL}, \cite{Ble}, \cite{BR}, \cite{BoR}, \cite{Ci}, \cite{Co},
\cite{Du}, \cite{Fe}, \cite{FV}, \cite{Galanis}, \cite{KMb},\cite{Os}, \cite{MZ}, \cite{MRM}-\cite{Mu10},  \cite{Pan},
\cite{Pau}, \cite{Ro9}-\cite{Rr},  \cite{Sly}, \cite{TL}, \cite{Tim}.

\section{Applications of the Potts model}

The most studied model of statistical mechanics is the Ising model (i.e. 2-state Potts model),
 in MathSciNet there are about 6000 papers devoted to the problems related to Ising model.
 In \cite{Mul} some physical motivations why the Ising model on a Cayley tree
  is interesting are given.
  In particular, this model plays a very special role in statistical mechanics and gives
  the simplest nontrivial example of a system undergoing phase transitions.

The \emph{Potts model} was introduced by
R.\,B.~Potts \cite{Potts} as a lattice system with $q\ge2$ spin
states and nearest-neighbor interaction, to generalize of
the Ising model.

The Potts model has been quickly picked up by a host of research in diverse areas.
Here we give brief (not complete\footnote{see https://worldwidescience.org/topicpages/q/q-state+potts+model.html
for a list and abstracts of works related to Potts model.}) review of
some interesting application of
the Potts model:

\subsection{Alloy behavior}

An alloy is\footnote{https://en.wikipedia.org/wiki/Alloy} a combination of metals or metals combined with one or more other elements.
Examples are combining the metallic elements gold and copper produces red gold, gold and silver becomes white gold, and silver combined with copper produces sterling silver.

By modeling of microstructural evolution one can study many metals processing
companies because the alloys are usually designed computationally by tailoring their
microstructural features (see \cite{AZ} and references therein).
These features entail grain size, particle/precipitate content,
recrystallization fraction, and crystallographic texture, among others \cite{Pol}.
Moreover, in many industrial thermo-mechanical processes, various annealing phenomena (for example:
recrystallization and grain growth) are incompletely understood.
 One of the computer simulation methods used to study such phenomena is the Monte Carlo
Potts model (MCPM). This model has been used to simulate annealing phenomena such as grain
growth in single- and two-phase polycrystalline materials, directional grain growth,
particle pinning, static and dynamic recrystallization, microstructure, abnormal and nanocrystalline
grain growth. The MCPM method has also demonstrated its applicability to modeling recrystallization
in aluminum alloys.

To simulate recrystallization in a particle containing alloy, the paper \cite{Rad} have coupled the finite element crystal plasticity
method with MCPM.

In \cite{AZ} a MCPM was used to model the primary recrystallization and grain
growth in cold rolled single-phase Al alloy.
The general vision of such modeling is to be able to optimize the
microstructural features (grain size, recrystallization fraction, and crystallographic texture)
during recrystallization and grain growth computationally in three dimensions in single-phase
Al alloys. This simulation provides beneficial tools for understanding annealing and related
phenomena in thermal treatments of rolled structures.

\subsection{Cell sorting}

In biology the cell \index{cell} (meaning ``small room") is the basic structural, functional, and biological
unit of all known organisms. A cell is the smallest unit of life.

Cell migration is a central process in the development and maintenance of multicellular organisms.
Tissue formation during embryonic development, wound healing and immune responses all
require the orchestrated movement of cells in particular directions to specific locations.
Cells often migrate in response to specific external signals, including chemical signals
and mechanical signals \cite{Mak}.

Cellular Potts models (CPMs) (see \cite{Vo} and references therein) are asynchronous
probabilistic cellular automata
developed specifically to model interacting cell populations. They are used
to the field of cell and tissue biology.  In particular, when
the details of intercellular interaction are essentially determined by the shape and the
size of the individual cells as well as the length of the contact area between neighboring
cells.

%
%
%
%
%

A CPM is a discrete-time Markov chain, where the
transition probabilities are specified with Gibbs measure corresponding to a Hamiltonian.

%
%

In \cite{AS} (see also the references therein) a CPM is considered,
to simulate single cell migration over flat substrates with variable stiffness.

In the migration the cells sense their surroundings and respond to different types of signals.
Cells preferentially crawl from soft to stiff substrates by reorganizing their
cytoskeleton from an isotropic to an anisotropic distribution of actin filaments.

The following configurations are studied \cite{AS}:
\begin{itemize}
\item[1.]\ \ a substrate including a soft and a stiff region,
\item[2.] \ \  a soft substrate including two parallel stiff stripes,
\item[3.] \ \  a substrate made of successive stripes with increasing stiffness
to create a gradient and
\item[4.] \ \ a stiff substrate with four embedded soft squares.
\end{itemize}
For each case it is evaluated the morphology of the cell, the distance covered,
the spreading area and the migration speed.

Cell sorting  is the process of taking cells from an organism
and separating them according to their type. The cells are labeled and tagged to identify areas of interest and their effect.
 They are separated based on differences in cell size, morphology (shape), and surface protein expression \cite{RKF}.


%
%
%

In the cell sorting while the surface-energy-driving mechanism is the same
as for grain growth, biological cells have generally a fixed
range of sizes. Thus the pattern cannot lose energy by
coarsening, since cells cannot disappear, similar constrained
evolution occurs in bubbles in magnetic films \cite{We}.

The differences in contact energies between
cells of different types (differential adhesion) cause cell
motion which reduces the pattern's energy.
To take this into account the authors of \cite{Graner} added
 an elastic-area constraint to the Hamiltonian of the Potts model:
introduce a symbol  $\tau$ for the cell type.
In the model there are three cell
types, ``light"$=l$, ``dark"$=d$ and ``medium"$=M$ that is $\tau \in \{l, d, M\}$. The
surface energy between two ce11s then depends on the
types of the cells. Each cell still has a unique spin $\sigma\in \Phi$,
and consists of all lattice sites with that
spin, but there may be many cells of each type, i.e., with
the same $\tau$.

The extended Potts model is  defined by the following
$$H(\sigma)=\sum_{\langle x, y\rangle}J(\tau(\sigma(x)), \tau(\sigma(y)))[1-\delta_{\sigma(x)\sigma(y)}]$$
$$+\lambda\sum_{{\rm spin \, types} \, \sigma}\left(a(\sigma)-A_{\tau(\sigma)}\right)^2 \theta(A_{\tau(\sigma)}),$$
where $\tau(\sigma)$ is the type associated with the cell $\sigma$ and
 $J(\tau, \tau')$ is the surface energy between types $\tau$ and
$\tau'$. $\lambda$ is a Lagrange multiplier specifying the strength of
the area constraint, $a(\sigma)$ the area of a cell $\sigma$, and $A_\tau$, the
target area for cells of type $\tau$.

Because of the surface energy,
each cell usually contains slightly fewer than $A_\tau$, lattice
sites. Moreover the biological aggregates are usually surrounded
by a fluid medium (i.e. $\tau =M$), e.g., culture solution, substrate,
or extracellular matrix, which is defined as a single
cell with associated type, interaction energies, and unconstrained
volume. Assume the target area $A_M$ of the medium
to be negative and to suppress the area constraint
include
$$\theta(x) = \left\{\begin{array}{ll}
0 \ \ \mbox{if} \ \ x < 0 \\[2mm]
1 \ \ \mbox{if} \ \ x > 0.
\end{array}\right.$$

It is interesting to determine whether a model of this type
exhibits biologically reasonable cell sorting. In coarsening,
patterns are usually characterized by their side and
area distributions and their moments, as well as the exponent
describing the average rate of area growth. In cell
sorting, the areas are approximately fixed so the area information
is not useful.

In \cite{Graner} the authors used Potts-model dynamics, with one Monte Carlo
time step defined to be 16 times the number of
spins in the array, but they suppress the nucleation of
heterogeneous spins by requiring that a lattice site flip
only to a spin belonging to one of its neighbors. It is found
a long-distance cell movement leading to sorting
with a logarithmic increase in the length scale of homogeneous clusters.
Sorted clusters then round. Moreover two successive phases are found:
a rapid boundary-driven creation of a low-cohesivity cell monolayer around
the aggregate, and a slower boundary-independent internal rearrangement.

\subsection{Financial engineering}

%

In \cite{Reichardt} based on a $q$-state Potts model a fast community detection algorithm is presented. Communities
in networks (groups of densely interconnected nodes that are only loosely connected to the rest of
the network) are found to coincide with the domains of equal spin value in the minima of a modified
Potts spin glass Hamiltonian. Comparing global and local minima of the Hamiltonian allowed for
the detection of overlapping (fuzzy) communities and quantifying the association of nodes to
multiple communities as well as the robustness of a community.

In \cite{Takaishi} a 3-state model based on the Potts model is proposed
to simulate financial markets. The three states are assigned to "buy", "sell" and "inactive" states. The model
shows the main stylized facts observed in the financial market: fat-tailed distributions
of returns and long time correlations in the absolute returns.

%

The work \cite{WW} uses the Potts model to simulate and characterize the time evolution of a market time series.
A two-dimensional 3-state Potts model is used to develop a stock price time series model.

This financial model imitates:
\begin{itemize}
\item[(i)] \ \ traders taking a selling position,
\item[(ii)] \ \ traders taking a buying position, and
\item[(iii)] \ \  traders taking no trading position,
\end{itemize}
 which are classified as type 1, type 2, and type 3, respectively.

It is assumed (in \cite{WW}) that stock price behavior is strongly affected
 by the number of traders $\omega^{(1)}(t)$ (traders of type 1),
 $\omega^{(2)}(t)$ (traders of type 2), and $\omega^{(3)}(t)$ (traders
of type 3).

Consider a single stock and assume that there are $L^2$ traders
of this stock who are located in a square lattice: $L\times L\subset \mathbb Z^2$.
Moreover, assume that each trader can trade a unit number of stock at each time $t\in \{1, 2, \dots, T \}$.
At each time $t$, the fluctuation of stock price process is strongly
influenced by the number of traders who take buying positions and the
number of traders who take selling positions. When the number of traders in selling positions is smaller than the number
of traders in buying positions, the stock price is considered low by market participants, and the stock price gradually
increases. The similar is true in the opposite case.
 In this proposed financial model, the clusters of parallel spins in the square-lattice Potts model are
designated groups of market traders acting together.



In \cite{XW} permutation entropy and sample entropy are
developed to the fractional cases, weighted fractional permutation entropy and
fractional sample entropy.  The effectiveness of these entropies as complexity
measures is analyzed by application to the logistic map, which is a typical
one-dimensional map creating chaos in some range of a parameter, such as the time series
of some stock market indices and the price evolution of an artificial stock market using the Potts model.  Moreover, the
numerical research on nonlinear complexity behaviors is compared between the returns
series of Potts financial model and the actual stock markets.

\subsection{Flocking birds}

Flocking\footnote{https://en.wikipedia.org/wiki/Flocking$_-$(behavior)}
is the behavior exhibited when a group of birds, called a flock,
are foraging or in flight. There are parallels with the shoaling behavior of fish,
the swarming behavior of insects, and herd behavior of land animals.

Examples: starlings are known for aggregating into huge flocks of hundreds to thousands of individuals,
murmur at ions, which when they take flight altogether, render large displays of intriguing swirling patterns in the skies above observers.

Mathematical models used to study the flocking behaviors of birds can also generally be applied to the ``flocking" behavior of other species.

In mathematical modeler, ``flocking" is the collective motion by a group of
self-propelled entities and is a collective animal behavior exhibited by
many living beings. The flocking transition is an out-of-equilibrium phenomenon
and abundant in nature (see \cite{CM} and references therein): from human
crowds, mammalian herds, bird flocks, fish schools to unicellular organisms such as amoebae,
bacteria, collective cell migration in dense tissues, and sub-cellular structures including cytoskeletal
filaments and molecular motors.

 The physics of flocking is also prevalent in nonliving substances such as
rods on a horizontal surface agitated vertically, self propelled
liquid droplets, liquid crystal hydrodynamics, and rolling colloids.

In basic models of flocking behavior the following three simple rules
are taken into account:
\begin{itemize}
\item \ \ Separation - avoid crowding neighbors (short range repulsion)
\item \ \ Alignment - steer towards average heading of neighbors
\item \ \ Cohesion - steer towards average position of neighbors (long range attraction)
\end{itemize}

With these rules, the flock moves in an extremely
realistic way, creating complex motion and interaction that
would be extremely hard to create otherwise.


In \cite{CM} the 4-state active Potts model (APM) is considered, which
has four internal states corresponding to four motion
directions and is defined on a two-dimensional lattice with
coordination number 4.
Its two main ingredients leading
to flocking are the local alignment interactions and self propulsion
via biased hopping to neighboring sites without
repulsive interactions.

%

 A local alignment
rule inspired by the ferromagnetic 4-state Potts model and self-propulsion via biased diffusion
according to the internal particle states leads to flocking at high densities and low noise.
The phase diagram of the APM is computed and  the flocking dynamics in the region is explored, in which
the high-density (liquid) phase coexists with the low-density (gas) phase and forms a fluctuating
band of coherently moving particles. Moreover, as a function of the particle self-propulsion velocity, a
novel reorientation transition of the phase-separated profiles from transversal to longitudinal band
motion is revealed, which is absent in the active Ising model.

\subsection{Flowing foams}

Foams (see \cite{CoA}) are complex fluids composed of gas bubbles
tightly packed in a surfactant solution.
In spite they generally consist only of Newtonian fluids,
foam flow obeys nonlinear laws. This can result from non-affine deformations of
the disordered bubble packing as well as from a coupling between the surface flow in
the surfactant mono-layers and the bulk liquid flow in the films, channels,
and nodes. A similar coupling governs the permeation of liquid through
the bubble packing that is observed when foams drain due to gravity.

In \cite{Sa1} the CPM successfully simulates drainage and shear in foams.
 The CPM is used to investigate instabilities due to the flow of a single large bubble in a dry,
 mono-disperse two-dimensional flowing foam.
 It is shown that as in experiments in a Hele-Shaw cell\footnote{The term Hele-Shaw cell is
 commonly used for cases in which a fluid is injected into the shallow geometry from above or below
 the geometry, and when the fluid is bounded by another liquid or gas.}, above a threshold velocity the large bubble moves
 faster than the mean flow. These simulations reproduce analytical and experimental predictions
 for the velocity threshold and the relative velocity of the large bubble,
 demonstrating the utility of the CPM in foam rheology studies.

\subsection{Image segmentation}

In digital imaging, a pixel (or picture element) is a physical point
in a raster\footnote{a raster is a dot matrix data structure
that represents a generally rectangular grid of colored points in a display medium.} image.
The pixel is the smallest controllable element of a picture represented on the screen.

More pixels typically provide more accurate representations of the original image.
In color imaging systems, a color is typically represented by three or four component
intensities such as red, green, and blue, or cyan, magenta, yellow, and black (which can be assigned to
states of the Potts model).

An image segmentation
is\footnote{https://en.wikipedia.org/wiki/Image$_-$segmentation} the process
of partitioning a digital image into multiple segments
(sets of pixels, also known as image objects).

The segmentation is needed to
simplify and/or change the representation of an image into something that
is more meaningful and easier to analyze. This is the process of assigning a label to every pixel in an image such that
 pixels with the same label share certain characteristics.

In image analysis the problems involving incomplete data arise when some part of the data
is missing or unobservable \cite{CF}. In this problems one wants to recover an original image
which is hidden.

%
%
%
%
%

In \cite{CF} a Markov model-based image segmentation is studied, which involves
hidden Markov random fields. The Potts model with external fields are used to better adequacy especially when the color
proportions are very unbalanced in the images to be recovered.

As usual the pixels set $S$ in image analysis is $\mathbb Z^2$ with the second order
neighborhood system: for each site, the neighbors are the eight sites surrounding it.

The image segmentation involves observed and unobserved  data to be recovered.

In \cite{CF} the unobserved data is modeled as a discrete Markov random field, where
  a commonly used distribution for the random field is the $q$-color Potts model with a non-zero
  external field. Moreover, an illustration of the advantages of introducing the external field is given with numerical experiments.

\subsection{Medicine}

Cancer  is a group of diseases characterized by the uncontrolled
growth of tumor\footnote{Uncontrolled cells may form a mass called a tumor}
cells that can occur anywhere
in the body.

 The paper \cite{KMM} contains a survey of mathematical models
that explicitly take into account the spatial architecture of three-dimensional
tumours and address tumour development, progression and response to treatments.
In particular, it discusses Potts models of epithelial acini, multicellular
spheroids, normal and tumour spheroids and organoids, and multi-component
tissues.

Here, following \cite{BJ}, \cite{LX}, \cite{Sun}, \cite{Sz}, we give a review of the application of
Potts model to cancer diseases.

In \cite{BJ} the authors discussed CPM's application in developmental biology, focusing on the
development of blood vessels, a process called vascular morphogenesis.
A range of models focusing on uncovering the basic mechanisms of
vascular morphogenesis are introduced: network formation and sprouting and then show how these
models are extended with models of intracellular regulation and with interactions
with the extracellular micro-environment. The integration of
models of vascular morphogenesis in several examples of organ development in
health and disease, including development, cancer, and age-related macular degeneration are reviewed.
The computational efficiency of the CPM and the
available strategies for the validation of CPM-based simulation models are presented.

Digital pathology imaging of tumor tissues, which captures histological details in high resolution, is
fast becoming a routine clinical procedure.

In \cite{LX} (see also references therein) the authors considered
 the problem of modeling a pathology image with irregular locations
of three different types of cells: lymphocyte, stromal, and tumor cells.
A novel Bayesian hierarchical model is given, which incorporates a hidden Potts model
to project the irregularly distributed cells
to a square lattice and a Markov random field
prior model  (Gibbs measures of the Potts model)  to identify regions in a heterogeneous
pathology image. The model allowed to quantify the interactions between different types of cells,
some of which are clinically meaningful.

The development of a tumor is initiated as the genomes of individual
cells in an organism become destabilized, which usually kills cells,
but in rare cases it modifies the properties
of the cell in a way that allows it to proliferate and form a
tumor.  Despite a growing wealth of available molecular data,
the growth of tumors, invasion of
tumors into healthy tissue, and response of tumors to therapies are still poorly understood \cite{Sz}.

In \cite{Sz} the authors review attempts to develop theoretical
frameworks for collective cell behavior during tumor development.
Mathematical descriptions of tumor growth and development range
from continuum-level descriptions of gene-regulatory
networks or tumor cell populations, to detailed, spatial models of
individual and collective cell behavior.
A cellular phenotype determines the success of a cancer cell in competition with
its neighbors, irrespective of the genetic mutations or physiological alterations that gave
rise to the altered phenotype. The CPM widely used to study the question what phenotypes
can make a cell ``successful" in an environment
of healthy and cancerous cells.   In particular, \cite{Sz} reviews use of the CPM for modeling
tumor growth, tumor invasion, and tumor progression.


In  \cite{Sun}, to simulate the spatiotemporal evolution of the tumor, a mathematical model is constructed.
Using both Potts model and nutrient competition a good visual simulation of tumor growth was provided,
reproducing experimental results which also shows that the tumor growth is sensitive to
the nutrient environment. It was noted that the results may have some medical signification:

(i)  Tumors grow exponentially in the beginning.
The tumor migrated toward the nutrient.

(ii) There are some differences in the dependence on nutrient between
malignant cells and healthy cells. It may allow to control the
nutrient environment in human host factitiously to avoid or cure cancer.

(iii) It may exist several critical nutrients which play a very
different even opposite role in the growth process of malignant cells.

\subsection{Neural network }

A neural network is a network or circuit of neurons,
or in a modern sense, an artificial neural network, composed of
artificial neurons or nodes. Thus a neural network is either a
biological neural network, made up of real biological neurons,
or an artificial neural network, for solving artificial intelligence problems.

%


The $q$-state Potts-glass model is a symmetric feedback
neural network with neurons of $q$ states \cite{Lo}, \cite{Xi}.
Each neuron can be modeled with a Potts spin $\sigma_i$, where subscript $i$
indicates neuron $i\in \{1, \dots, q\}$, that may represent
the color or the shade of grey of a pixel in a pattern.

In \cite{Aq} the critical properties of the Potts model with $q=3$ and 8 states in one-dimension
on directed small-world networks are investigated. The Potts model on these networks presents
 a second-order phase transition with a new set of critical exponents for $q=3$.
 For $q=8$ the system exhibits only a first-order phase transition.

The authors of \cite{Lo} have studied the $q$-states Potts generalization
of the Hopfield neural network evolving in parallel. Precise results are obtained on the
asymptotic number of stored independent and identically distributed patterns so that either the
patterns are fixed points of the dynamics or corrupted patterns get attracted to the original
ones after one or several steps of the dynamics.

In \cite{Xi} the authors have studied the $q$-state Potts-glass neural network with the pseudo-inverse rule.
It is found that there
exists a critical point of $q_c=14$, below which the storage capacity and the retrieval quality can be
greatly improved by introducing the pseudo-inverse rule. It was shown that the dynamics of the neural networks constructed
with the two learning rules respectively are quite different; but however, regardless of the learning rules, in the
$q$-state Potts-glass neural networks with $q\geq 3$ there is a common novel dynamical phase in which the spurious
memories are completely suppressed. This property has never been noticed in the symmetric feedback neural
networks. Free from the spurious memories implies that the multi-state Potts-glass neural networks would not
be trapped in the meta-stable states, which is a favorable property for their applications.

In very recent paper \cite{TL} using the techniques of neural networks, the three - dimensional,
5-state ferromagnetic Potts model on the cubic lattice and two-dimensional 3-state
anti-ferromagnetic Potts model on the square lattice are studied.
The whole or part of the theoretical
ground state configurations of the studied models are considered as the training sets.

The results of the three-dimensional Potts model
imply that the neural networks approach is as efficient as the traditional method since the signal
of a first order phase transition, namely tunneling between two channels, determined by the neural networks
method is as strong as that calculated with the Monte Carlo technique. Furthermore, the outcomes
associated with the considered two-dimensional Potts model indicate even little partial information of the ground
states can lead to conclusive results regarding the studied phase transition. These results
demonstrate that the performance of neural networks, using the
theoretical ground state configurations as the training sets, is impressive.

\subsection{Phases transitions}

In fact, the Potts model firstly appeared in statistical physics \cite{Ba}, \cite{Wu}.
In previous sections of this paper we gave strong mathematical presentation of this
model for physical point of view.

It is known\footnote{https://en.wikipedia.org/wiki/Potts$_-$model} that as a model
of a physical system, the Potts model is simple, but useful as
a model system for the study of phase transitions.
For example, the ferromagnetic Potts model on $\mathbb Z^2$ has a first order transition if $q > 4$.
In case $q = 4$ a continuous transition is observed, as in the Ising model where $q = 2$.
Other use of Potts model is found through the model's relation to percolation problems
and the Tutte and chromatic polynomials found in combinatorics \cite{Sa}, \cite{Wu}.

For integer values of $q\geq 3$, the model displays the phenomenon of 'interfacial adsorption'
with intriguing critical wetting properties when fixing opposite boundaries in two different states.

The $q$-state Potts model with $q\ge3$ and an external field $\alpha\in\mathbb{R}$ on Cayley tree was
considered in \cite{Peruggi1}, \cite{Peruggi2}) using physical argumentation. In particular,
they identified the critical temperature point $T_{\rm c}$
and also suggested an explicit critical boundary in the phase
diagram for $\alpha\ge0$.

It should be stressed that the phase transition occurring at these
critical boundaries is not of type ``uniqueness/non-uniqueness'',
with which we are concerned in the present book, but in
fact the so-called ``order/disorder'' phase transition. The latter
was studied rigorously in \cite{Galanis} in connection with the computational complexity of
approximating the partition function of the Potts model. The useful
classification of critical points deployed in \cite{Galanis} is
based on the notion of \emph{dominant phase}; in particular, a
critical point is determined, from this point of view
as a threshold beyond which only ordered phases are dominant.

\subsection{Political trends}
In \cite{NO} the Potts model is applied to Twitter data \index{Twitter data} related to political trends.
Twitter is a micro blogging environment where users post small messages,
or Twitts, \index{Twitts}  depicting their likes and dislike
towards a certain topic, e.g. candidates to the next political elections.

On the base of several electoral events and assuming a stationary regime,
the authors of \cite{NO} found the following:

(i) the maximization of the entropy singles out a microscopic
model (single-Twitt-level) that coincides with a $q$-state Potts model having suitable
couplings and external fields to be determined via an inverse problem from the two
sets of data;

(ii) correlations decay as $1/N_{e}$, where $N_e$ is a small fraction of the mean
number of Twitts;

(iii) the simplest statistical models that reproduce these correlations
are the multinomial distribution, characterized by $q$ external fields,
and the mean field
Potts model, characterized by one coupling;

(iv) this coupling
turns out to be always close to its critical value.

\subsection{Protein family}

A protein is a polypeptide chain consists of a sequence
of amino acids. The sequence of these amino acid units
and one additional state for gaps or empty spaces. The
gap state must be available for an amino acid to move
when they fold into three dimensional structures to form
domains.

A protein family\footnote{https://en.wikipedia.org/wiki/Protein$_-$family} is a
group of evolutionary-related proteins.
Usually a protein family has a corresponding gene family, in which each
gene encodes a corresponding protein with a 1:1 relationship.
The evolutionary history of a protein family is typically represented by
a phylogenetic tree\footnote{https://en.wikipedia.org/wiki/Phylogenetic$_-$tree}.

Proteins in a family descend from a common ancestor and typically
have similar three-dimensional structures, functions, and significant sequence similarity.
Two segments of a DNA may have shared ancestry because of three phenomena: either a speciation event,
or a duplication event, or else a horizontal gene transfer event. Sequence homology
is the biological homology between DNA, RNA, or protein sequences,
defined in terms of shared ancestry in the evolutionary history of life.

 Potts model is used to describe the sequence variability of sets of evolutionarily
  related protein sequences  (i.e. homologous
protein families). The notion of homologous protein families implies
that present sequences derive from a common ancestor (see \cite{RE} and references therein).

%

In \cite{RE} the authors gave a principled way to correct for phylogenetic effects in the inference
of Potts models from sequence data. Although the standard technique to account for these effects in
co-evolutionary analysis relies on an empirical re-weighting of sequences. The method of \cite{RE}
aims at doing so using the phylogenetic tree as well as an evolutionary model. The global nature of Potts models
implies that the evolutionary model used should depend on the full sequence. An inference
scheme is given that takes the phylogeny of a protein family into account in order to correct biases in
estimating the frequencies of amino acids. Using artificial data, it is shown that a Potts model inferred
using these corrected frequencies performs better in predicting contacts and fitness effect of mutations.

\subsection{Protein folding}

Protein folding is the physical process meaning that
a protein chain acquires its native three-dimensional structure,
a conformation that is usually biologically functional, in an expeditious and reproducible manner.

This subsection is based on \cite{SS}.

A Hamiltonian representing interacting
amino acids on a lattice is given by
$$H(\sigma)=-{1\over 2} \sum_{i\ne j}J_{ij}\delta_{\sigma_i,\sigma_j}-\sum_i h_i \delta_{\sigma_i,\sigma_{q_0}},$$
where the external field $h_i$ favors variables to align in $\sigma_{q_0}$; that is forming
domains of amino acid $\sigma_{q_0}$.
Here $J_{ij}$ is the pairwise interaction strength or the exchange parameter between
variables $\sigma_i$ and $\sigma_j$ at two different lattice
sites $i$ and $j$. Note that the model allows each lattice site to have one out of $q$ different states.

The Kronecker delta function and
the negative sign in front of the first term of the Hamiltonian favor to have
same amino acids at sites $i$ and $j$.
Due to this pairwise
attractive interaction between same residues, domains of
same amino acids is expected to be formed by folding
the amino acid chain into a three dimensional structure.
By the $q$-state Potts model, this structural formation is characterized by ``ordering"
of spin variables $\sigma$ in
a specific state.

\begin{figure}
\includegraphics[width=7cm]{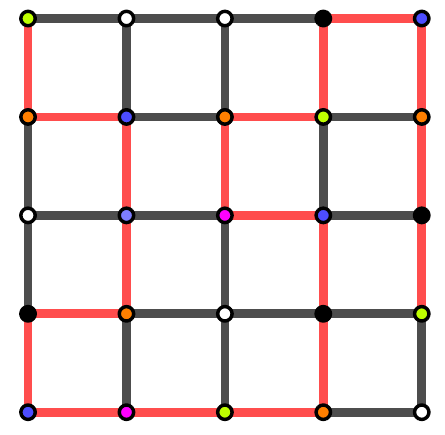}
\caption{ The two-dimensional protein is
represented as a chain of amino acids (colored dots)
connected by peptide bonds (red lines). The different
colored dots occupied at lattice sites represent the
different amino acids, the white dots represent the
empty sites on the two dimensional lattice.}\label{pf}
\end{figure}

In case when the
entropy dominates at higher temperatures, the chain
stretches inside the lattice (see Fig. \ref{pf}) by moving amino acids to
empty sites. When the energy dominates at lower
temperatures, favored by the Potts model, same color
amino acids cluster together by moving them closer.
This clustering and stretching of same colors due to the
competition between the energy and the entropy
represents the protein folding and un-folding.

The general investigation of
protein folding problem requires to solve the inverse Potts
model, i.e. finding the coupling constants $J_{ij}$ and $h_i$
from the data base.

The main protein folding problem is understanding the question of how a protein's amino acid sequence
dictates its structure. The various interaction parameters inside the protein and the local environment are
directly responsible for the folding.  In \cite{SS} the authors studied the
thermodynamic properties of the protein folding (the Potts model)
using a statistical mechanics
approach. Namely, converting the interacting amino acids into an effectively non-interacting model using a
mean-field theory, the Helmholtz free energy  is found. Then by investigating the Helmholtz free energy,
the properties of protein folding transition is qualitatively studied. It is shown that the protein folding
phase transition is a strongly first order.

Potts models of protein sequence co-variation are statistical models constructed from the
pair correlations observed in a multiple sequence alignment of a protein family.
In (\cite{LH} see also references therein) the Potts models are reviewed to predict protein
structure and sequence-dependent conformational free energy landscapes, to survey protein fitness
landscapes and to explore the effects of epistasis on fitness.

\subsection{Sociology}

In \cite{BEP},  \cite{Sc}, \cite{Schulze} an application of Potts model to studies of human behavior is given:
Nobel laureate Thomas Schelling published a seminal paper \cite{Sc} which in addition to organized and economic
explanations considers the possibility of micro-motive explanations for racial segregation. \index{racial segregation}
The premise is that individual decisions to avoid minority status (or to require being in a minority of some
minimum size) could lead to the macro-effect of segregation. Schelling places vacancies, stars, and zeros randomly on a
checkerboard and then iteratively considers the happiness of the stars and zeros with their local neighborhoods, moving
an unhappy star or zero to the nearest vacant spot that meets their happiness criteria.
Schelling works through many different experiments to come up with some very
compelling results on segregation. From this model it seems that people do consciously
or subconsciously segregate themselves from people who are different than they are.

The work \cite{MO} gives
a similar premise to Schellings (that micro-motive explanations can lead to immigrant ghettos) with a more Potts-like
model where the Hamiltonian measures the happiness of individuals with their neighbors, the temperature is viewed as
a social temperature where warmer temperatures reflect facilitation of integration and assimilation, and at each step in the
simulation two neighbors are able to exchange places with a probability based on the likelihood of the new state with respect
to the current state. The paper \cite{Schulze} extends the work of \cite{MO} to address up to seven different ethnic groups.

Following \cite{BL} one can construct a Potts model for simulating human behavior in the following
way. Use a lattice to depict your neighborhood, city, business, or any other venue
in which people interact with one another. To have more groups, we can have elderly people,
college roommates, families with teenagers, and
families with small children. Assume that members of each of these groups of people
are living together in a brand new development. Label the elderly with a 1, the
college roommates with a 2, the families with teenagers with a 3, and the families with
small children with a 4.
The members of these groups have preferences about who they live near. For
example, the elderly do not want to live next to the college roommates because of the
large parties that they tend to throw. The couples with small children might want to live
next to one another so that their kids can play together without going far from home.
The Hamiltonian for this experiment would measure overall happiness as opposed
to energy. Outside forces might be the price of other houses in other neighborhoods,
proximity to work, or how much people like their current house. The Metropolis
Algorithm could then be run to develop higher probabilities for lattice states with higher
overall happiness. Eventually, we would likely see preferences playing out in the form of
segregation.
This is just a rough sketch of a Potts model scenario. This may give an appreciation
for the versatility of the Potts model when it comes to real world situations.

\subsection{Spin glasses}

 A glassy system is complex system that exhibits a very slow dynamics
 that prevents it from reaching the equilibrium state. Examples are
 real glasses, spin glasses, supercooled liquids, polymers, granular material,
colloids, ionic conductors, orientational glasses, and vortex
glasses \cite{FC}.

 A spin glass is a model\footnote{https://en.wikipedia.org/wiki/Spin$_-$glass} of
 a certain type of magnet. Magnetic spins are the orientation of
 the north and south magnetic poles in three-dimensional space. In ferromagnetic solids,
 component atoms magnetic spins all align in the same direction.
 Spin glasses are contrasted with ferromagnetism as ``disordered" magnets in which their atoms
  spins are not aligned in a regular pattern.

  The complex internal structures that arise within spin glasses are termed ``metastable"
     because they are ``stuck" in stable configurations other than the lowest-energy configuration.
     The mathematical complexity of these structures is difficult but fruitful to study experimentally or in simulations; with applications to physics, chemistry, materials science and artificial neural networks in computer science.

The Potts model of spin glass is defined by the Hamiltonian \cite{FC}:
 $$H(\sigma, s)=-qJ\sum_{\langle i, j\rangle}\left[\delta_{\sigma_i\sigma_j}(\epsilon_{ij}s_is_j+1)-2\right],$$
 where associated with each lattice site is an Ising spin $s_i\in \{-1,1\}$
and a $q$-state Potts spin $\sigma_i\in \{1,\dots,q\}$. The sum is
extended over all nearest-neighbor sites, $\epsilon_{i,j}=\pm 1$ is a
random quenched variable, and $J$ is the strength of interaction.

The model is a superposition of a ferromagnetic $q$-state
Potts model and a $\pm J$ Ising spin glass model \cite{MP}.

In \cite{FC} it is shown that this model exhibits for all $q$ a spin glass
transition at $T_{SG}(q)$ and a percolation transition at higher temperature $T_P(q)$. It is shown that for all values
of $q>1$ at $T_P(q)$ there is a thermodynamic transition in the universality class of a ferromagnetic $q$-state Potts
model. Moreover, the efficiency of the cluster dynamics is compared with that of standard spin-flip dynamics.

There are variations of the Potts spin glass model, which has been studied extensively in
the physics literature (see, e.g., \cite{Cal}, \cite{MP}, \cite{Pan} and references therein).


\subsection{Storage capacity}

 Storage capacity measures how much data a computer
 system may contain. For an example, a computer with a 500GB hard drive has a
 storage capacity of 500 gigabytes.

 A Potts unit (spin value)   can be regarded (see \cite{N}) in the neuroscience
 context as representing a local subnetwork or cortical patch of real neurons,
 endowed with its set of dynamical attractors, which span
different directions in activity space, and are therefore converted to the states of the
Potts unit (which is defined precisely as having states pointing each along a different
dimension of a simplex). One can define the model as an auto-associative network of Potts
units interacting through tensor connections. The memories are stored in the weight
matrix of the network and they are fixed: each
memory $\mu$ is a vector or list of the states taken in the overall activity configuration by
each unit.

 Take each Potts unit to have $S$ possible active states, labeled e.g. by
the index $k$, as well as one quiescent state, $k = 0$, when the unit does not participate in
the activity configuration of the memory. Using this model in \cite{N} the storage capacity of
the Potts network is studied.  The storage capacity calculation, performed using
replica tools, is limited to fully connected networks, for which a Hamiltonian
can be defined. To extend the results to the case of intermediate partial
connectivity,  the self-consistent signal-to-noise analysis is derived for the
Potts network; and the implications for semantic memory in
humans is discussed.

In \cite{Kan} the theory of neural networks is extended to include discrete neurons  with $q$, $q\geq 2$ discrete states.
The dynamics of such systems are studied by using Potts model. The maximum number of storage patterns
is found to be proportional to $Nq(q-1)$, where $q$ is the number of Potts states and $N$ is the size of the network.

It is known (see \cite{Kr}) that the capacity to store information in any device, and in particular the capacity to store
concepts in the human brain, is limited. In \cite{Kr} it is shown in a minimal model of semantic
memory, and in progressive steps, how one can expect the storage capacity to behave
depending on the parameters of the system.
It was deduced the minimum requirements of any model of this kind in
order to have a high capacity. The calculation specifies that in the Potts model the number of concepts that can be stored
is neither linear nor an arbitrary power  of the number $S$ of values a feature can
take, but quadratic in $S$.

\subsection{Symmetric channels}

Symmetric channels on $q$ symbols have the state space $\{1,\dots,q\}$ and $q\times q$-matrices: \index{symmetric channels}
\begin{equation}\label{me8}
{\bf M} =\left(\begin{array}{ccccc}
1-(q-1)\delta& \delta&\delta&\dots &\delta\\[2mm]
\delta& 1-(q-1)\delta&\delta&\dots&\delta\\[2mm]
\dots&\vdots &\dots & \vdots&\dots\\[2mm]
\delta&\delta&\delta&\dots&1-(q-1)\delta
\end{array}\right),
\end{equation}
with $\lambda_2(M)=1-q\delta$.

Depending on the sign of $\lambda_2(M)$ we distinguish between ferromagnetic Potts models
where $\lambda_2(M)>0$, and anti-ferromagnetic models where $\lambda_2(M)<0$.
In case $1-(q-1)\delta=0$ we obtain the model of proper colorings of the tree:

\begin{equation}\label{me9}
{\bf M} =\left(\begin{array}{ccccc}
0& (q-1)^{-1}&(q-1)^{-1}&\dots &(q-1)^{-1}\\[2mm]
(q-1)^{-1}& 0&(q-1)^{-1}&\dots&(q-1)^{-1}\\[2mm]
\dots&\vdots &\dots & \vdots&\dots\\[2mm]
(q-1)^{-1}&(q-1)^{-1}&(q-1)^{-1}&\dots&0
\end{array}\right).
\end{equation}

Several bounds for the reconstruction problem (for $q$-state Potts models on a Cayley tree of order $k\geq 2$)
are known (see \cite{Mos1}-\cite{Mos}, \cite{Sly}):
\begin{itemize}
\item \ \ If $k\lambda_2^2(M) > 1$ then the reconstruction problem is solvable.

\item \ \ If $k|\lambda_2(M)|\leq 1$, then the reconstruction problem is unsolvable.

\item \ \ If $k\lambda_2(M) > 1$ and $q$ is sufficiently large, then
the reconstruction problem is solvable.

\item \ \
If $k{(1-q\delta)^2\over 1-(q-2)\delta}\leq 1$ then the reconstruction problem is unsolvable.
\end{itemize}
As it was mentioned above the most general result on reconstruction is the Kesten-Stigum bound \cite{Ke}
which says that reconstruction holds when $k\lambda_2^2 >1$.

In \cite{Sly} the author proved the first exact reconstruction
threshold in a non-binary model establishing the Kesten-Stigum bound for the
$3$-state Potts model on regular trees of large degree. The Kesten-Stigum
bound is not tight for the $q$-state Potts model when $q\geq 5$.
Moreover,  asymptotics for these reconstruction thresholds are determined.

It it was studied in genetics (see e.g. \cite{Cav}, \cite{St1}) an information flow can be used to represent
propagation of a genetic property from ancestor to its
descendants.  In communication theory,
this process represents a communication network on the tree where information is
transmitted from the root of the tree. Moreover, the process was studied in statistical
physics, see e.g. \cite{Bl}, \cite{YH}. More precisely, in statistical physics a
non-uniqueness of the Gibbs measure means that for all $n$ there exists $\sigma_n$ such
that the distribution of $\sigma_{\rho}$ given $\sigma_n$ has total
variation distance at least $\delta > 0$ from uniform. This is a weaker condition
than reconstruction solvability and it was studied in statistical physics
for Ising and Potts models (see \cite{Ge}, \cite{Ha}, \cite{P}).

The crucial role of the reconstruction problem in Phylogeny
was demonstrated in \cite{Mos4} and \cite{Mos5}.

\subsection{Technological processes}

Technology means activities, because it is in constant development.
In general technological processes can be design processes,
making processes and processes in the phase of using and assessing technology \cite{Vr}.

Ostwald ripening\footnote{https://en.wikipedia.org/wiki/Ostwald$_-$ripening} is a phenomenon
 observed in solid solutions or liquid sols that describes the change of an inhomogeneous
 structure over time, i.e., small crystals or sol particles dissolve, and redeposit onto
  larger crystals or sol particles. The Ostwald ripening is generally found in water-in-oil
  emulsions, while flocculation is found in oil-in-water emulsions.

Some analytic methods used to model processes such as
grain growth and Ostwald ripening by many simplifying
assumptions.

In \cite{Sr} the $q$-state Potts model was adapted to
study grain growth (development). By this model it was treated the next level of complexity
by incorporating grain boundary topology. Moreover, the growth kinetics of
a two dimensional, connected assembly of mutually interacting grains are described.

Now following  \cite{Ti} we give one more application of the Potts model  \cite{Ti}.
Markets often demand small quantities of many different components with highly specialized
performance requirements. Therefore, the development cost
of each component cannot be amortized over large numbers
and becomes prohibitive. Good predictive tools are necessary for reduction of the development cost.
It have been proved that the Potts model is a powerful and useful tool for simulating a
wide variety of microstructural evolution problems both during
processing and during component use. In \cite{Ti} the use
and application of the Potts model to Ostwald ripening is
presented.

Here, to see how the Potts model appears,  we give phase equilibria characteristics of \cite{Ti}:

Consider $\mathbb Z^2$ as the lattice of the $q$-state Potts model.
 A canonical ensemble of sites, which may be visualized
as building blocks of the microstructure, populates a simulation
lattice. These sites possess certain energies based on their inherent
characteristics and on their interactions with their neighboring
domains and evolve to minimize the total free energy of
the system.

\begin{figure}[th]
\includegraphics
[height=8cm]{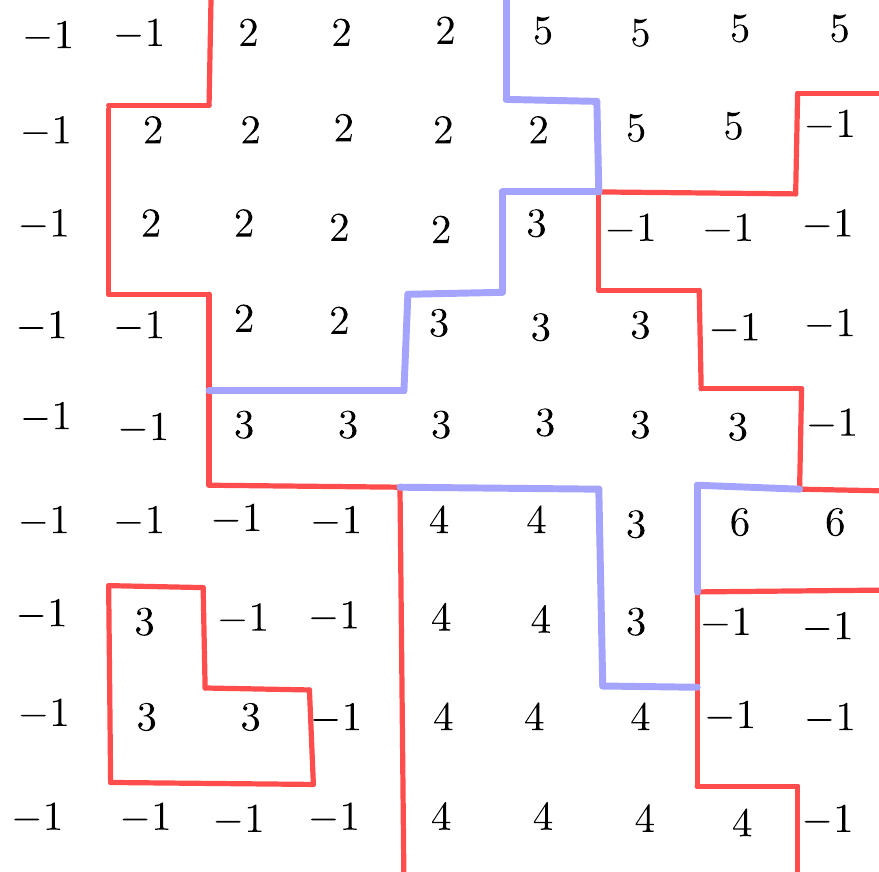}%
\caption{Illustration of the technique used to digitize
microstructures. The grains in a liquid with grain boundaries (blue) and
solid-liquid interfaces (red) are shown. The $B$-sites are all the sites labeled
``-1" and $A$-sites are labeled with positive integers.}\label{fgr}
\end{figure}

The evolution mechanism is that each domain can
be modified based on the energetics of the modification being
considered. Then a digitized microstructures can be represented on  $\mathbb Z^2$
with a periodic boundary conditions. The components necessary
for Ostwald ripening can be incorporated into the Potts
model by populating the simulation lattice with a canonical
ensemble having two components designated as $A$ and $B$. The
$A$-sites are the primary building blocks of the grains and the
$B$-sites of the matrix, as shown in Fig. \ref{fgr}.

The two components,
$A$ and $B$, must separate into two phases, grains dispersed in a
matrix. The two phases were generated by defining bond energies
between the $A$- and $B$-sites, so that the components would segregate into two phases.
The $A-B$ bonds are assigned
higher energies than $A-A$ or $B-B$ bonds. Thus the system minimizes
its energy by segregating the $A$- and $B$-sites into two
phases.

\subsection{Wetting transition}

Here we follow \cite{LM} (see also references therein).
 The study of liquid behavior on a surface is a topic that is
receiving increasing interest due to its diverse applications.

Examples of application are: microfluidics for biotechnology, textile lubrication, and
self-cleaning. The inspiration arose in nature, more precisely
from the plant (the lotus), whose leaves
have a natural hydrophobicity.

Based on the surface structure of the leaf,
similar artificial rough surfaces were produced. A drop of water on
this type of surface can exhibit two states of wettability:
\begin{itemize}
\item Wenzel state, in which the drop of water penetrates
the surface cavities;
\item Cassie-Baxter (CB) state, in which the drop remains on the top of the pillars,
as in the lotus leaf.
\end{itemize}

The Wenzel approach considers that the
liquid fills up completely the grooves on the rough surface
increasing thus the total area of the liquid-solid interface.

The CB approach assumes that air pockets are trapped inside the grooves, and
thus, that the liquid only contacts the solid at the top of
the asperities. Thus, the drop sits on a composite surface
comprised of solid and trapped air.

The energy barrier between the CB-Wenzel state transition is given as
a product of the following factors: the pillar
height, pillar thickness, and the area that the droplet
occupies between the pillars. This energy barrier separating the
states may be large enough to prevent a spontaneous transition from occurring.

The work \cite{LM} is based on a two-dimensional Potts model of
a droplet placed on a regularly patterned surface of pillars
of different heights and interspacing.

The model is defined as follows.
Consider a square lattice, where to each
network site, $i$, a value $\sigma_i\in \Phi=\{1,\dots, q\}$ is given.
The set $\Phi$ is labels associated with the
different media: air, liquid, solid etc.

The Hamiltonian (the total energy) of the
system is
$$H(\sigma)={1\over 2}\sum_{\langle i,j\rangle}\gamma_{\sigma_i\sigma_j}\left(1-\delta_{\sigma_i\sigma_j}\right)
+\lambda\left(A_{\sigma=1}-A_T\right)^2+mg\sum_i\delta_{1\sigma_i}h_i,$$
where $\sigma=\{\sigma_i: \, i\in \mathbb Z^2\}$ the configuration (the labels on the sites of the lattice) and
$\gamma_{\sigma_i\sigma_j}$ is the interfacial energy between neighboring sites belonging to different
media \cite{Graner}. The second term corresponds to an elastic compression energy: $A_{\sigma=1}$
 is the drop area, $A_T$ is a target area,
and $\lambda$ is a Lagrange multiplier related to the inverse of the
liquid compressibility \cite{Id}. This term is required to avoid
the droplet disappearance by energy minimization. The last
term is the contribution from the gravitational energy, where
$g$ is the acceleration of gravity\footnote{The gravity can be
neglected in case of very small droplets.} and $h_i$ is the vertical
coordinate of a label belonging to the droplet defined as having
unit mass $(m_{\sigma=1} = 1)$.

In \cite{Sha} (see also the references therein) the
transition and the influence of the energy barrier on the
transition through thermodynamics and molecular dynamics models.

In \cite{LM} Potts model simulation compared with the studies developed \cite{Sha}, \cite{Shah} employing
thermodynamics and two-dimensional molecular dynamics.

\section*{ Acknowledgements}

I thank my coauthors who collaborated with me
to develop the theory of Gibbs measures: L.V.Bogachev, G.I.Botirov,
M.Cassandro, D.Gandolfo, N.N.Ganikhodjaev,  F.Haydarov, F.Henning,
O.Khakimov, R.Khakimov,  S.Kissel, C.K\"ulske,
A.LeNy, L.Liao,  I.Merola, F.M.Mukhamedov,  M.Rahmatullaev,
J.Ruiz, Yu.M.Suhov, J.P.Tian, Y.Velenik and many others.

\end{document}